\documentclass{rspublic_mod}
\usepackage{graphicx}

\usepackage[T1]{fontenc}
\usepackage[latin9]{inputenc}
\usepackage{amssymb}
\usepackage[T1]{fontenc}
\usepackage[latin9]{inputenc}
\usepackage{amstext}

\begin{document}

\title{The Dirichlet Series for the Liouville Function and the Riemann Hypothesis}

\author [K. Eswaran] {K. Eswaran}

\affiliation{Srinidhi Institute of Science and Technology, Hyderabad, India.\footnote{The penultimate  sub-section `An Intuitive Analogy to Understand the Formal Results' and the Appendix I, II have been introduced here to reach out to a larger audience. In this Version, a paragraph on application of the iterated logarithm has been added. Typos have been removed and the Appendix 5 is slightly modified and a new Appendix 6 has been added.   Correspondence may be addressed to: keswaran@sreenidhi.edu.in OR kumar.e@gmail.com } }

\label{firstpage}

\maketitle
\begin{abstract}{Liouville function, multiplicative function, factorization into primes, Mobius function, Riemann hypothesis}
\noindent This paper investigates the analytic properties of the Liouville
function's Dirichlet series that obtains from the function $F(s)\equiv $ $%
\zeta (2s)/\zeta (s)$, where $s$ is a complex variable and $\zeta (s)$\ is
the Riemann zeta function. The paper employs a novel method of summing the
series by casting it as an infinite number of sums over sub-series that
exhibit a certain symmetry and rapid convergence. In this procedure, which
heavily invokes the prime factorization theorem, each sub-series has the
property that it oscillates in a predictable fashion, rendering the analytic
properties of the Dirichlet series determinable. With this method, the paper
demonstrates that, for every integer with an even number of primes in its
factorization, there is another integer that has an odd number of primes
(multiplicity counted) in its factorization. Furthermore, by showing that a sufficient condition
derived by Littlewood (1912) is satisfied, the paper demonstrates that the
function $F(s)$ is analytic over the two half-planes $Re(s)>1/2$ and $%
Re(s)<1/2$. This establishes that the nontrivial zeros of the Riemann
zeta function can only occur on the critical line $Re(s)=1/2$.

\end{abstract}
\section{Introduction}
This paper investigates the behaviour of the Liouville function which is
related to Riemann's zeta function, $\zeta (s)$, defined by%
\begin{equation}
\zeta (s)=\sum\limits_{n=1}^{\infty }\frac{1}{n^{s}},
\end{equation}

\noindent where $n$ is a positive integer and $s$ is a complex number,
with the series being convergent for $Re(s)>1$. This function has
zeros (referred to as the trivial zeros) at the negative even integers $%
-2,-4,\ldots $. It has been shown\footnote{%
This was first proved by Hardy (1914).} that there are an infinite number of zeros on the line at $Re(s)=1/2 $. Riemann's Hypothesis (R.H.) claims that these are all the nontrivial zeros of the zeta function. The R.H. has eluded proof to date, and this paper demonstrates that it is resolvable by tackling the Liouville
function's Dirichlet series  generated by $F(s)\equiv \zeta (2s)/\zeta (s)$%
, which is readily rendered in the form%
\begin{equation}
F(s)=\sum\limits_{n=1}^{\infty }\frac{\lambda (n)}{n^{s}},
\end{equation}

\noindent where $\lambda (n)$ is the Liouville function defined by $\lambda
(n)=(-1)^{\Omega (n)}$, with $\Omega (n)$ being the total number of prime
numbers in the factorization of $n$, including the multiplicity of the
primes. We would also need the summatory function $L(N)$, which is defined  as the partial sum up to $N$ terms of the following series:
\begin{equation}
L(N)=\sum\limits_{n=1}^{N}\,\lambda (n) \tag{1.2b}
\end{equation}
 
Since the function $F(s)$ will exhibit poles at the zeros of $\zeta
(s)$, we seek to identify where $\zeta (s)$ can have zeros by examining the
region over which $F(s)$ is analytic. By demonstrating that a sufficient
condition, derived by Littlewood (1912), for the R.H. to be true is
indeed satisfied, we show that all the nontrivial zeros of the
zeta function occur on the `critical line' $Re(s)=1/2$.

Briefly, our method consists in judiciously partitioning the set of positive
integers (except 1) into infinite subsets and couching the infinite sum in (1.2) into sums over
these subsets  with each resulting sub-series being uniformly convergent. This method of considering a slowly  converging series as a sum of many sub-series was previously used by the author in problems where Neumann series were involved Eswaran (1990)).

In this paper  we break up the sum of the Liouville function into sums over many sub-series whose behaviour is predictable. It so turns out that one prime number $p$ (and its powers)  which is associated  with a particular sub-series controls the behaviour of that sub-series.

Each sub-series is in the form of rectangular functions (waves) of unit amplitude but ever increasing periodicity and widths - we call these `harmonics' - so that every prime number is thus associated with such harmonic rectangular functions which then play a role in contributing to the value of $L(N)$. It so turns out that if N goes from N to N+1, the new value of $L(N+1)$ depends solely on the factorization of N+1, and the particular harmonic that contributes to the change in $L(N)$ is completely determined by this factorization. Since prime factorizations of numbers are uncorrelated, we deduce that the statistical distribution of $L(N)$ when N is large is like that of the cumulative sum of N coin tosses, (a head contributing +1 and a tail contributing -1), and thus logically lead to the final conclusion of this paper.

We found a new method of factoring every integer and placing it in an exclusive subset, where it and its other members form an increasing sequence which in turn  factorize alternately into odd and even factors; this method exploited the inherent symmetries of the problem and was very useful in the present context. Once this symmetry was recognized, we saw that it was natural to invoke it in the manner in which the sum in (1.2) was performed. We may view the sum as one over subsets of series that exhibit convergence even outside the domain of the half-plane $Re(s)>1$. We  were rewarded, for following the procedure pursued in this paper, with the revelation that the Liouville function (and therefore the zeta function) is controlled by innumerable rectangular harmonic functions whose form and content are now precisely known and  each of which is associated with a prime number and all prime numbers play their due role. And in fact all harmonic functions associated with prime numbers below or equal to a  particular value N determine L(N).

When we are oblivious to the underlying symmetry being alluded to here, we render the summation in (1.2) less tractable than necessary. This is precisely what happens when we perform the
summation in the usual manner, setting $n=1,2,3,...$ in sequence.

In addition to establishing that the evenness, denoted as $(+1)$, and oddness, denoted as $(-1)$, of the number of prime factors of \emph {consecutive positive integers} behave like the results from the tossing of an ideal coin, we also establish that the sequence of $+1$'s and $-1$'s can never be cyclic (see Appendix III). In Appendix IV, we offer an intuitive proof of the claim that the sequence
of $\lambda $'s occurring in $L(N)$ for large $N$ behaves like coin tosses.
This is followed in Appendix V by a formal, rigorous arithmetic proof of the same,
along with a determination of the asymptotic behavior of $|L(N)|$---thus
completing the validation of the Riemann Hypothesis. As a final confirmation, by using Kolmogrov's law of the iterated algorithm, we show (in the end of Section 5), that the `width' of the Critical Line, as expected, vanishes to zero. 

The main paper and the Appendices I to V we concern ourselves with the mathematical proof of the R.H. In Appendix VI, we perform a numerical analysis and provide supporting empirical evidence that is consistent with the
formal theorems that were key to establishing the correctness of the RH. By performing this exhaustive numerical analysis and statistical study we obtain a clearer understanding of the Riemann Problem and its resolution.

\section{Partitioning the Positive Integers into Sets}
\setcounter{equation}{2}

The Liouville function $\lambda (n)$ is  defined over the set of positive integers $n$ as  $\lambda(n) = (-1)^{\Omega(n)}$, where $\Omega(n)$ is the number of prime factors of $n$, multiplicities included. Thus $\lambda(n)=1$
when $n$  has an even number of prime factors and $\lambda(n)=-1$
when it has an odd number of prime factors. We define $\lambda(1)=1$. It is a completely arithmetical function obeying $\lambda(mn)=\lambda(m)\lambda(n)$ for any two positive integers $m \,,\,n$.

 We shall consider subsets of positive integers such as $%
\{n_{1},n_{2},n_{3},n_{4},...\}$ arranged in increasing order and are such that
their values of $\lambda $ alternate in sign: 
\begin{equation}
\lambda (n_{1})=-\lambda (n_{2})=\lambda (n_{3})=-\lambda (n_{4})=...
\end{equation}

\noindent It turns out that we can label such subsets with a triad of
integers, which we now proceed to do. To construct such a labeling scheme,
consider an example of an integer $n$ that can be uniquely factored into primes as
follows:%
\begin{equation}
n=p_{1}^{e_{1}}p_{2}^{e_{2}}p_{3}^{e_{3}}...p_{L}^{e_{L}}p_{i}p_{j}
\end{equation}

\noindent where $p_{1}<p_{2}<p_{3}...<p_{L}<p_{i}<p_{j}$ are prime numbers
and the $e_{k},k\in \{1,2,3,...,L\}$ are the integer exponents of the
respective primes, and $p_{L}$ is the largest prime with exponent exceeding $%
1$, the primes appearing after $p_{L}$ will have an exponent of only one and there may a finite number of them, though only two are shown above. Integers of this sort, with at least one multiple prime factor are
referred to here as Class I integers. In contrast, we shall refer to
integers with no multiple prime factors as Class II integers. A typical
integer, $q$, of Class II may be written%
\begin{equation}
q=p_{1}p_{2}p_{3}...p_{j}p_{L},
\end{equation}

\noindent where, once again, the prime factors are written in increasing
order.

We now show how we construct a labeling scheme for integer sets that
exhibit the property in (2.3) of alternating signs in their corresponding $%
\lambda $'s. First consider Class I integers. With reference to (2.4), we
define integers $m,p,u$ as follows:%
\begin{equation}
m=p_{1}^{e_{1}}p_{2}^{e_{2}}p_{3}^{e_{3}}...p_{L-1}^{e_{L-1}};\quad
p=p_{L};\quad u=p_{i}p_{j}.
\end{equation}

\noindent In (2.6), $m$ is the product of all primes less than $p_{L}$,the 
largest multiple prime in the factorization, and $u$ is the product of
all prime numbers larger than $p_{L}$ in the factorization. Thus the Class I
integer $n$ can be written%
\begin{equation}
n=mp^{e_{L}}u
\end{equation}%
Hence we will label this integer $n$ as $(m,p^{e_{L}},u)$,using the triad of numbers%
$(m,p,u)$ and the exponent $ e_{L}$. It is to be noted that $u$ will consist of prime factors
all larger than $p$, and $u$ cannot be divided by the square of a prime number.

Consider the infinite set of integers, $P_{m;p;u}$, defined by%
\begin{equation}
P_{m;p;u}=\{mp^{2}u,mp^{3}u,mp^{4}u,...\}
\end{equation}

\noindent The Class I integer $n$ necessarily belongs to the above set
because $e_{L}\geq 2$. Since the consecutive integer members of this set
have been obtained by multiplying by $p$, thereby increasing the number of
primes by one, this set satisfies property (2.3) of alternating signs of the
corresponding $\lambda $'s. Note that the lowest integer of this set $%
P_{m;p;u}$ of Class I integers is $mp^{2}u$.

We may similarly form a series for Class II integers. The integer $q$ in (2.5)
may be written $q=mpu$, with $m=p_{1}p_{2}p_{3}...p_{j}$, $p=p_{L}$, and $u=1$%
. This Class II integer is put into the set $P_{m;p;u}$ defined by%
\begin{equation}
P_{m;p;1}=\{mp,mp^{2},mp^{3},mp^{4},...\}.
\end{equation}

\noindent The set containing Class II integers is distinguished by the facts
that $u=1$ for all of them, their largest prime factor is always $p$ and none
 of them can be divided by the square $\pi^2$ of a prime number 
$\pi$ such that  $ \pi  < p $; in other words the factor $m$  cannot
 be divided by the square of a prime. In this set, too, the $\lambda $'s
alternate in sign as we move through it and so property (2.3) is satisfied.
Again, note that the lowest integer of this set $P_{m;p;1}$ is the Class II
integer $mp$, all the others being Class I.

 In what follows, we shall find it handy to refer to the set of
ascending integers comprising $P_{m;p;u}$  as a `tower'. It is
important to distinguish between a tower (or set) described by a triad like 
$(m,p,u)$  and an integer belonging to that set. It is worth repeating that the set or tower of
Class I integers described by the label $(m,p,u)$  is the infinite
sequence $\{mp^{2}u,mp^{3}u,mp^{4}u,...\}$, the first element of
which is  $mp^{2}u$  and all other members of which are $mp^{k}u$,
 where $k>2$. A set or tower containing a Class II integer
described by $(m,p,u=1)$ is the infinite sequence 
$\{mp,mp^{2},mp^{3},...\}$, Eq.(2.9), of which only the first element $mp$ %
 is a Class II integer and all other members, $mp^{k}$,
where $k\geq 2$, are Class I, because the latter have exponents greater than 1.  For convenient reference, we shall
refer to the first member of a tower as the base integer or the base of the
tower. It is also worth noting that when we refer to a triad like $%
(m,p^{k},u)$, where $k>1$, we are invariably referring to
the integer $mp^{k}u$  and not to any set or tower. Labels for sets
do not contain exponents; only those for integers do. Of course, the particular integer
 $(m,p^{k},u)$ belongs to the set or tower $(m,p,u)$.

Two simple examples illustrate the construction of the sets denoted by $%
P_{m;p;u}$:

\noindent Ex. 1: The integer $2160$, which factorizes as $2^{4}\times
3^{3}\times 5$, is clearly a Class I integer since it is divisible by the
square of a prime number---in fact there are two such numbers, $2$ and $3$%
---but we identify $p$ with $3$ as it is the larger prime. It is a member of
the set $P_{16;3;5}=\{16\times 3^{2}\times 5,16\times 3^{3}\times 5,16\times
3^{4}\times 5,16\times 3^{5}\times 5,...\}$.

\noindent Ex. 2: The integer $663$, which factorizes as $3\times 13\times 17$%
, is a Class II integer because it is not divisible by the square of a prime
number. It belongs to the set $P_{39;17;1}=\{39\times 17,39\times
17^{2},39\times 17^{3},....\}$.

Note that two different integers cannot share the same triad.%
\footnote{%
The integer represented by the triad $(m,p^{r},u),$ is the product $mp^{r}u$%
, which obviously cannot take on two distinct values.} And two different
triads cannot represent the same integer.\footnote{%
Suppose two different triads $(m,p^{r},u)$ and $(\mu ,\pi ^{\rho },\nu )$
represent the same integer, say $n$. Then we must have $mp^{r}u=\mu \pi
^{\rho }\nu =n.$It follows that at least two numbers of the tetrad $%
\{m,p,r,u\}$ must differ from their counterparts in the tetrad $\{\mu ,\pi
,\rho ,\nu \}$. Since the factorization of $n$ is unique, this is impossible.%
} Thus the mapping from a triad to an integer is one-one and onto. A formal proof is in the Appendix.

The following properties of the sets $P_{m;p;u}$ may be noted:

(a) The factorization of an integer $n$ immediately determines whether it is
a Class I or a Class II type of integer.

(b) The factorization of integer $n$ also identifies the set $P_{m;p;u}$ to
which $n$ is assigned.

(c) The procedure defines all the other integers that belong to the same set
as a given integer.

(d) Every integer belongs to some set $P_{m;p;u}$ (allowing for the
possibility that $u=1$)\ and only to one set. This ensures that,
collectively, the infinite number of sets of the form $P_{m;p;u}$ exactly
reproduce the set of positive integers $\{1,2,3,4,....\}$, without omissions
or duplications.

Our procedure, taking its cue from the deep connection between the
zeta function and prime numbers, has constructed a labeling scheme that
relies on the unique factorisation of integers into primes. In what follows,
we shall recast the summation in (1.2) into one over the sets $P_{m;p;u}$.%
 The advantage of breaking up the infinite sum over all positive
integers into sums over the $P_{m;p;u}$ sets will soon
become clear. 

\section{An Alternative Summation of the Liouville Function's Dirichlet
Series}
\setcounter{equation}{9}

We shall now implement the above partitioning of the set of all positive
integers to examine the analytic properties of $F(s)$ in (1.2). We shall
rewrite the sum in (1.2) into an infinite number of sums of sub-series,
ensuring that each sub-series is uniformly convergent even as $s\rightarrow
0$.

We begin, however, by assuming that $Re(s)>1$, which makes the series
in (1.2) absolutely convergent. We write the right hand side in sufficient
detail so that the implementation of the partitioning scheme becomes
self-evident:%
\begin{eqnarray}
F(s) &=&1+\sum\limits_{r=1}^{\infty }\frac{\lambda (2^{r})}{2^{rs}}%
+\sum\limits_{r=1}^{\infty }\frac{\lambda (3^{r})}{3^{rs}}%
+\sum\limits_{r=1}^{\infty }\frac{\lambda (5^{r})}{5^{rs}}%
+\sum\limits_{r=1}^{\infty }\frac{\lambda (2\times 3^{r})}{2^{s}3^{rs}} 
\notag \\
&&+\sum\limits_{r=1}^{\infty }\frac{\lambda (7^{r})}{7^{rs}}%
+\sum\limits_{r=1}^{\infty }\frac{\lambda (2\times 5^{r})}{2^{s}5^{rs}}%
+\sum\limits_{r=1}^{\infty }\frac{\lambda (11^{r})}{11^{rs}}%
+\sum\limits_{k=2}^{\infty }\frac{\lambda (2^{k}\times 3)}{2^{ks}3^{s}} 
\notag \\
&&+\sum\limits_{r=1}^{\infty }\frac{\lambda (13^{r})}{13^{rs}}%
+\sum\limits_{r=1}^{\infty }\frac{\lambda (2\times 7^{r})}{2^{s}7^{rs}}%
+\sum\limits_{r=1}^{\infty }\frac{\lambda (3\times 5^{r})}{3^{s}5^{rs}}%
+\sum\limits_{r=1}^{\infty }\frac{\lambda (17^{r})}{17^{rs}}  \notag \\
&&+\sum\limits_{r=1}^{\infty }\frac{\lambda (19^{r})}{19^{rs}}%
+\sum\limits_{k=2}^{\infty }\frac{\lambda (2^{k}\times 5)}{2^{ks}5^{s}}%
+\sum\limits_{r=1}^{\infty }\frac{\lambda (3\times 7^{r})}{3^{s}7^{rs}}%
+\sum\limits_{r=1}^{\infty }\frac{\lambda (2\times 11^{r})}{2^{s}11^{rs}} 
\notag \\
&&+\sum\limits_{r=1}^{\infty }\frac{\lambda (23^{r})}{23^{rs}}%
+\sum\limits_{r=1}^{\infty }\frac{\lambda (2\times 13^{r})}{2^{s}13^{rs}}%
+\sum\limits_{k=2}^{\infty }\frac{\lambda (2^{k}\times 7)}{2^{ks}7^{s}}%
+\sum\limits_{r=1}^{\infty }\frac{\lambda (29^{r})}{29^{rs}}\notag \\
&&+\sum\limits_{r=1}^{\infty }\frac{\lambda (2\times 3\times 5^{r})}{2^{s}3^{s}5^{rs}}+ \cdots %
\end{eqnarray}

We have explicitly written out a sufficient number of terms of the right
hand side of (1.2) so that those corresponding to each of the first 30
integers are clearly visible as a term is included in one (and only one) of
the sub-series sums in (3.10). On the right hand side, the second term contains
the integers $2,4,8,16....$; the third contains $3,9,27,...$; the fourth
contains $5,25,125,...$; the fifth contains $6,18,54,...$; sixth contains $%
7,49,...$; the seventh contains $10,50,...$; the eighth contains $11,121,...$%
; the ninth contains $12,24,48,...$; and so on. Note that in the ninth,
fifteenth, and twentieth terms the running index is deliberately switched
from $r$ to $k$ to alert the reader to the fact that the summation starts
from $2$ and not from $1$ as in all the other sums. (Note that, in the ninth
term, the Class I integer $n=12=2^{2}\times 3$ is assigned to the set $%
P_{1;2;3}=\{2^{2}\times 3,2^{3}\times 3,2^{4}\times 3,...\}$ and not to the
set $P_{4;3;1}=\{2^{2}\times 3,2^{2}\times 3^{2},2^{2}\times 3^{3},.\}$,
because the first term identifies $p$ as $2$ and $u$ as $3$ where as the second  term onwards $3$ has exponents, which violates our rules of precedence and would be an illegitimate
assignment given our partitioning rules.)

The sub-series in (3.10) have one of two general forms:%
\begin{equation}
\sum_{r=1}^{\infty}\:\frac{\lambda(m.p^{r})}{m^{s}.p^{rs}}\:\:=\frac{\lambda(m.p)}{m^{s}.p^{s}}\:[1-\frac{1}{p^{s}}+\frac{1}{p^{2s}}-\frac{1}{p^{3s}}+\cdots+\frac{(-1)^{X}}{p^{Xs}}+\cdots ] \notag
\end{equation}
or
\begin{equation}
\sum_{k=2}^{\infty}\:\frac{\lambda(m.p^{k}.u)}{m^{s}.p^{k}.u^{s}}\:\:=\frac{\lambda(m.p^{2}.u)}{m^{s}.p^{2s}.u^{s}}\:[1-\frac{1}{p^{s}}+\frac{1}{p^{2s}}-\frac{1}{p^{3s}}+ \cdots +\frac{(-1)^{X}}{p^{Xs}} +\cdots ]
\end{equation}
The above geometric series occurring within square brackets in the above two equations can actually be summed (because they are convergent),(see Whittaker and Watson) but we will refrain from doing so,
and (1.2) can be rewritten as%

\begin{equation}
F(s)=\sum\limits_{m}\sum\limits_{p}\sum\limits_{u} \,F^{T}_{m;p;u}(s)+\sum\limits_{m}\sum%
\limits_{p}\, F^{T}_{m;p;1}(s),
\end{equation}

\noindent where the first group of summations pertains to Class I integers $n$
characterized by the triad $(m,p^{k},u),(k \ge 2)$ and the second group
pertain to those integers which are characterized by set $(m,p^{k},1), (k \ge
1) $ the first member in the set is a Class II integer and others Class I.

In the above  we have defined the function $F^{T}_{m;p;u}(s)$  of the complex variable $s$ which is a sub-series involving terms over only the tower $(m,p,u)$   for a Class
I integer as follows%
\begin{equation}
F^{T}_{m;p;u}(s)= \sum\limits_{k=2}^{\infty} \frac{\lambda (mp^{k}u)}{m^{s}p^{ks}u^{s}},
\end{equation}

\noindent and the function $F^{T}_{m;p;1}(s)$ of the complex variable $s$ which is a sub-series involving terms over only the tower $(m,p,1)$   whose 1st term is a Class II integer as%
\begin{equation}
F^{T}_{m;p;1}(s)=\sum\limits_{r=1}^{\infty }\frac{\lambda (mp^{r})}{m^{s}p^{rs}}
\end{equation}

With the understanding that when $u=1$ we use the function in (3.14) instead
of (3.13), we may write $F(s)$ as%
\begin{equation}
F(s)=\sum\limits_{m}\sum\limits_{p}\sum\limits_{u}\,F^{T}_{m;p;u}(s).
\end{equation}

Comparing the above Eq.(3.15) with Eq(3.10) one can easily see that each  term which appears as a summation in (3.10) is actually a sub-series over some tower which we denote as $F^{T}_{m;p;u}(s) $ in (3.15). So we see that $F(s)$  has been broken up into a number of sub-series. The important point to note is that the $\lambda$ value of each term in the sub-series changes its sign from +1 to -1 and then back to +1 and -1 alternatively. Therefore if the starting value of $\lambda$ at the base was +1 then the cumulative contribution of this tower (sub series)  to $L(N)$ as N, the upper bound, increases from $N$ to $N+1,\,N+2,\,N+3,\,.... $ will fluctuate between  be 0 and 1. For some other tower whose base value of $\lambda$ is $-1$ its cumulative contribution to L(N) will fluctuate between 0 and $-1$;  these cumulative contributions can be represented in the form of a rectangular wave as shown in Figure 1. 

We have arrived at a critical point in our paper. We have cast the original
function $F(s)\equiv \zeta (2s)/\zeta (s)$ as a sum of functions of $s$.
Since the triad $(m,p^{k},u)$ uniquely characterises all integers, the
summations over $m,p,k$ and $u$ above are equivalent to a summation over all
positive integers $n$, as in (1.2), though not in the order $n=1,2,3,4,...$.
The manner in which the triads were defined ensures that there are neither
any missing integers nor integers that are duplicated.(See Theorems A and B in Appendix II.)

Although we did not explicitly do it, we mentioned in passing
that the sum over $k$ in (3.13)and (3.14) is readily performed since it is a
geometric series (see (3.11)) that rapidly converges. This is true not merely for 
$Re(s)>1$  but also as  $Re(s)\rightarrow 0 $. Whether 
 $F(s)$ converges when the summation is carried out over all the
towers $(m,p,u)$ and, if so, over what domain of $s$ is the central
question that we seek to answer in the next section. The answer to which as we shall see determines the analyticity of F(s) and thus resolves the Riemannian Hypothesis.

We can recast (3.15), still in the domain $Re(s)>1$, in the form%
\begin{equation}
F(s)=\sum\limits_{n=1}^{\infty }\frac{g(n)}{n^{s}},
\end{equation}

\noindent where $g(n)$ is a function appropriately defined below.

By construction, every $n$ in the above summation can be written as%
\begin{equation}
n=\mu \pi ^{\rho }\nu ,
\end{equation}

\noindent where $\mu $, $\pi $, and $\rho $ are positive integers, $\pi $ is
the largest prime in the factorization of $n$, with either (\textit{i}) an
exponent $\rho \geq 2$, and $\nu $ is the product of primes larger than $\pi 
$ but with exponents equal to $1$ (for Class I integers) or (\textit{ii}) it
is the largest prime factor with $\rho =1$ and $\nu =1$ (for Class II
integers).

We define $g(n)$ as follows:

\begin{eqnarray}
&=&\lambda (mp^{k}u)\quad \text{if }\mu =m\text{ and }\pi =p\text{ and }\nu
=u\neq 1\text{ and }\rho =k > 1 \quad \quad (3.18a) \notag \\
g(n) &=&\lambda (mp^{k})\quad \text{if }\nu =u=1\text{ and }\pi =p\text{\
and }\rho =k \ge 1 \quad\quad\quad  \quad \quad \quad \quad  (3.18b) \notag \\
&=&0\quad \text{otherwise.} \quad \quad\quad\quad\quad\quad\quad%
\quad \quad \quad \quad \quad \quad \quad \quad \quad \quad \quad \quad \quad \quad \quad (3.18c) \notag
\end{eqnarray}

The factors $m^{s}p^{ks}u^{s}$ and $m^{s}p^{ks}$ in the denominators of (3.13)
and (3.14) are simply $n^{s}$, where $n$ is the integer characterized by the $%
(m,p^{k},u)$ triad (with $u=1$ in the latter case).

\section{Calculation of the summatory Liouville function L(N)}

We are now in a position to examine the summatory Liouville function $L(N)$ by actually summing up the individual contribution from each sub-series.

To do all this systematically, we will explicitly illustrate the process starting from $N=1,2,3...$ up to $N=15$. Each of these numbers is factored and expressed uniquely as a triad.
The N=1 is a constant term, which is the trivial $(1,1,1)$,  then the next number $N=2 = (1,2,1)$, is contained in the tower shown below the one corresponding to $N=1$; and $N=3 =(1,3,1)$, is the tower below the previous; $4 =(1,2^2,1)$ however $4$ is already contained in the tower $(1,2,1)$ as its second member; the next N's: $5,6,7$, give rise to the new towers $(1,5,1), (2,3,1),(1,7,1); 8$ of course is the third member of the old tower $(1,2,1)$ similarly  $9$ is the 2nd member of $(1,3,1)$. After this the new towers which make their appearance are: $10=(2,5,1), 11=(1,11,1),13=(1,13,1),14=(2,7,1)$ and $15=(3,5,1)$. Figure 1 shows these and numbers up to N=30.
\begin{figure}
\includegraphics[width=13cm,height=13cm]{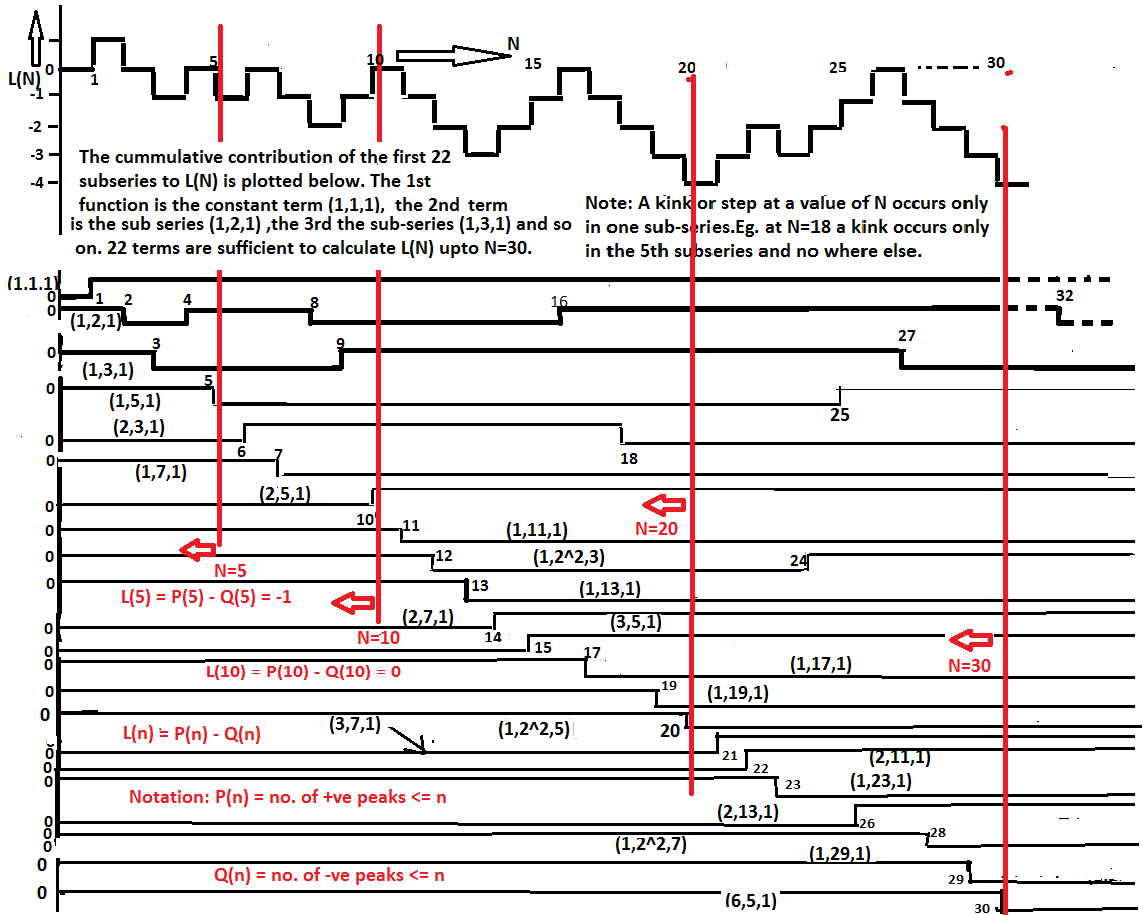}
\caption{The cumulative sum, L(N) (see top), is obtained by `filling up' slots in various towers from the bottom up until we have exhausted all N integers. }
\label{fig:1}       
\end{figure}
Now each  tower $(m,p,u)$ contributes to $L(N)$ (consider N fixed in the following) according to the following rules:

(i) A particular tower will contribute only if its base number is less 
 than or equal to N, i.e. $m.p.u \le N$ 
 
(ii) And the contribution $C$ to $L(N)$ from this particular tower will be exactly as follows:

Case A; Class II integer $(u=1)$

 $C=\Sigma_{r=1}^{R} \lambda(m.p^{r}.1)$ , where R is the largest integer such that $m.p^{R} \le N$ 
 
 Case B; Class I integer $(u>1)$
 
 $C=\Sigma_{k=2}^{K} \lambda(m.p^{k}.u)$, Where K is the largest integer such that $m.p^{K}.u \le N$ 
 
Now since each successive $\lambda $ changes sign  from $+1$ to $-1$ or vice a versa, the contributions of each tower can be thought of as a rectangular wave of 
ever-increasing width but constant amplitude $-1$ or $+1$, see Figure 1. 
 
 To find the value of L(N), (N fixed), all we need to do is count the jumps of each wave: as we move from N=0 a jump upwards is called a positive peak, a jump downwards is a negative peak. Draw a vertical line at N, we are assured that it will hit one and only one peak (positive or negative) in one of the sub-series; then count the total number of positive peaks P(N) and negative peaks Q(N), of the waves on and to the left of this vertical line, then $L(N) = P(N) - Q(N)$; the reason for this rule will be clear after the next section.

 For an example, take $N=5$. There is a positive peak for the constant term (1,1,1), the next wave (1,2,1) contributes one negative peak (at 2) and a positive peak (at 4), the wave (1,3,1) contributes a $-1$ peak (at 3) and (5,1) contributes a $-1$ peak (at 5). Thus a total of three negative peaks and two positive peaks add up to give $L(5)= -1$, which is of course correct. Now if we take $N=10$, and draw a vertical line at N=10, looking at this line and to its left we see that there are  additionally three positive and two negative peaks thus adding this contribution of +1 to the previously calculated value $ L(5)$ we get $L(10)=0$. (Two red vertical lines just just beyond N=5 and N=10 are drawn for convenience.) Now if we wish to compute L(15) we see that there are three more negative peaks and two positive peaks thus giving a value $L(15)=-1$. Counting the peaks further on it is easy to check that $L(N)$ is correctly predicted for every value of N up to 30 and in particular, $L(20)=-4$, $L(26)=0$ and $L(30)=-4$.

In summary, to calculate L(N) we merely need to count the negative and positive peaks of the waves on N and to the left of N. In the figure we have drawn a number of waves and labeled the tower to which each belongs using a triad of numbers. They are sufficient for one to easily calculate L(N) up to N=30 and check them out by comparing the numbers with the plot of L(N) shown on the top of the figure.

We turn to a more fundamental point: We show, in Section 6, that, for sufficiently large N (see Appendix IV),  the distribution of the value of  L(N) is equivalent to that obtained from summing the distribution of N coin tosses.

\setcounter{equation}{18}
\section{Analyticity of $F(s)$ and the Riemann Hypothesis}

We now utilize a technique introduced by Littlewood (1912) to examine the analyticity of the function $F(s)$.  In this, we follow
the treatment of Edwards (1974, pp. 260-261). The series in (3.16) can be
expressed as the integral%
\begin{equation}
F(s)=\int\limits_{0}^{\infty }x^{-s}dG(x), (Re(s) > 1),  \tag{5.19}
\end{equation}

\noindent where $G(x)=\int\nolimits_{0}^{x}dG$ is a step function that is
zero at $x=0$ and is constant except at the positive integers, with a jump
of $g(n)$ at $n$. The value of $G(n)$ at the discontinuity, at an integer $n$%
, is defined as $(1/2)[G(n-\epsilon )+G(n+\epsilon )]$, which is equal to $%
\sum\nolimits_{j=1}^{n-1}g(j)+(1/2)g(n)$. Assuming $Re(s)>1$,
integration by parts yields%
\begin{eqnarray}
F(s) &=&\int\limits_{0}^{\infty }d[x^{-s}G(x)]-\int\limits_{0}^{\infty
}G(x)d[x^{-s}]  \quad\quad \quad\quad (5.20)\notag\\
&=&\underset{X\rightarrow \infty }{\lim }\quad \lbrack
X^{-s}G(X)+s\int\limits_{0}^{X}G(x)x^{-s-1}dx  \notag \\
&=&s\int\limits_{0}^{\infty }G(x)x^{-s-1}dx,  \quad\quad(5.21) \notag
\end{eqnarray}

\noindent where the last step follows from the fact that $\mid G(X)\mid \le X $,
which implies that $X^{-s}G(X)\rightarrow 0$ as $X\rightarrow \infty $. We
further observe, following Littlewood (1912), that as long as $G(X)$ grows
less rapidly than $X^{a}$ for some $a>0$, the integrals in (5.21) and in the line 
preceding it converge for all $s$ in the half-plane $Re%
(a-s)<0$, that is, for $Re(s)>a$. By analytic continuation, $F\left(
s\right) $ converges in this half-plane. Since this result will be important
in what follows, we record it here.

\ 

\emph{Theorem 1 [Littlewood (1912)]: When }$G(X)$\emph{\ grows less rapidly
with }$X$\emph{\ than }$X^{a}$\emph{\ for some }$a>0$\emph{, }$F(s)$\emph{\
is analytic in the half-plane }$Re(s)>a$\emph{.}

\ 

We shall now demonstrate that the sufficient condition stated in Theorem 1
is satisfied for a specific value of $a$ that settles the Riemann
Hypothesis.  

From now on we revert to the original definitions of the sequence $ g(n) \equiv \lambda(n)$ and $G(N) \equiv L(N) $ as defined in Eq. (1.2) but write them in the forms derived in Section 3. Hence our  definition of $G(N)$ becomes %
\begin{equation}
G(N)=\sum\limits_{n=1}^{N}g(n),  \tag{5.22}
\end{equation}

\noindent and we may rewrite $G(N)$ as%
\begin{equation}
G(N)=\sum\limits_{m}\sum\limits_{p}\sum\limits_{u}\sum\limits_{k}\left[
(1-\delta _{u,1}).(1-\delta _{k,1})\,\lambda (mp^{k}u)+\delta _{u,1}\lambda (mp^{k})\right] , 
\tag{5.23}
\end{equation}

\noindent where $\delta _{u,1}$ and $\delta _{k,1}$ are Kronecker deltas (e.g. $\delta _{u,1}=1$ if $%
u=1$ and $0$ otherwise). 
The summations over $m$, $p$, $k$, and $u$ in (5.23) are undertaken with the
understanding that the triads $(m,p^{k},u)$ will only include integers $%
n\leq N$. 
Since the summation over $k$ is over an individual tower(if we keep (m,p,u) fixed we can write(5.23) as 
\begin{equation}
G(N)=\sum\limits_{m}\sum\limits_{p}\sum\limits_{u}\,F^T_{m,p,u}(s=0) , 
\tag{5.23b}
\end{equation}
This is nothing but Eq.(3.15) evaluated from each subseries $F^T_{m,p,u}(s)$ by making $s \rightarrow 0$. 

Of course, what we have called $G(N)$ is really the summatory
Liouville function, $L(N)$, defined earlier by (1.2b):

\begin{equation}
L(N)=\sum\limits_{n=1}^{N}\lambda (n).  \tag{5.24}
\end{equation}

Expression (5.23) is crucial because, in the light of Theorem 1, its behaviour
will determine the validity of the Riemann Hypothesis. Every term in the
summation in (5.23) is either $+1$ or $-1$. We need to determine, for given $N$,
 how many terms contribute $+1$ and how many $-1$, and then determine how
the sum $G(N)$ varies with $N$.

As we go through the list $n=1,2,3,\cdots, N$, we are assigning the integers to
various sets of the kind $P_{m;p;u}$. To use our terminology of towers, we
shall be `filling up' slots in various towers from the bottom up until we
have exhausted all $N$ integers. (When $N$ increases, in general, we shall
not only be filling up more slots in existing towers but also adding new
towers that were previously not included.) So the behaviour of $G(N)$ is
determined by how many of the numbers that do not exceed $N$
contribute $+1$ and how many $-1$.

It is convenient to identify the $\lambda $\ of an integer by the triad
which uniquely defines that integer. To avoid abuse of notation, we shall
denote the value $\lambda (n)$ in terms of the $\lambda$-value of the base
integer of the tower to which $n$ belongs. We will define the $\lambda $ of
the base of a tower in uppercase, as $\Lambda (m,p,u)$. In other words if $%
n=(m,p^{\rho },u)$ then it will belong to a tower whose base number is $%
n_{B}\equiv(m,p^{\kappa },u)$, where $\kappa =2\,\,if$ $\,u\neq 1\,$and $\,\kappa
=1\,if\,$\ $u=1.$ Now we define $\Lambda (m,p,u)=\lambda (n_{B})=\lambda
(mp^{\kappa }u)=\lambda (m)\lambda (p^{\kappa })\lambda (u)$, since the $%
\lambda $ of a product of integers is the product of the $\lambda $ of the
individual integers. Of course, once we know $\lambda (n_{B})$ we will know
the $\lambda $ of all other numbers belonging to the tower because they
alternate in sign.

To determine the behaviour of $G(N)$, the following theorem is important.

\emph{Theorem 2: For every integer that is the base integer of a tower
labeled by the triad }$(m,p,u)$\emph{, and therefore belonging to the set }$%
P_{m;p;u}$, \emph{there is another unique tower labeled by the triad }$%
(m^{\prime },p,u)$ \emph{and therefore belonging to the set }$P_{m^{\prime
},p,u}$\emph{\ \ with a base integer for which }$\Lambda (m^{\prime
},p,u)=-\Lambda (m,p,u)$.

\underline{Proof}:

Let us write the integers at the base of a tower in the form $n=mp^{\rho }u$
described by the triad $(m,p,u)$, where we shall assume that $\rho =2$ if $%
u\neq 1$ and $\rho =1$ if $u=1$. These correspond to the smallest members of
sets of Class I and Class II integers, respectively, which are the integers
of concern here. In the constructions below, we shall multiply (or divide) $%
m $ by the integer $2$. Since $2$ is the lowest prime number, such a
procedure does not affect either the value of $p$ or $u$ in an integer and
so we can hold these fixed.

We begin by excluding, for now, triads of the form $(1,p,1)$, integers which
are single prime numbers. We allow for this in Case 3 below.

\underline{Case 1}: Suppose $m$ is odd. We choose $m^{\prime }=2m$, then 

$\Lambda (m^{\prime },p,u)=\Lambda (2m,p,u)=-\Lambda(m,p,u)$. We may say that
 $(m,p,u)$ and $(m^{\prime },p,u)$ are `twin' pairs
in the sense that their $\Lambda $s are of opposite sign. Note that $(m,p,u)$
and $(m^{\prime },p,u)$ are integers at the base of two different towers;
they are not members of the same tower. (Recall that the members of a given
tower are constructed by repeated multiplication with $p$.)

\underline{Case 2}: Suppose $m$ is even. In this case, we need to ascertain
the highest power of $2$ that divides $m$. If $m$ is divisible by $2$ but
not by $2^{2}$, assign $m^{\prime }=m/2$. (So $m=6$ gets assigned to $%
m^{\prime }=3$, and $m=3$, by Case 1 above, gets assigned to $m^{\prime }=6$%
.) More generally, suppose the even $m$ is divisible by $2^{k}$ but not by $%
2^{k+1}$, where $k$ is an integer. Then, if $k$ is even, assign $m^{\prime
}=2m$; and if $k$ is odd, assign $m^{\prime }=m/2$. (So $m=12=2^{2}\times 3$
gets assigned to $m^{\prime }=2^{3}\times 3=24$. And, in reverse, $%
m=24=2^{3}\times 3$ gets assigned to $m^{\prime }=24/2=12$.)

Thus for odd $m$ the following sequence of pairs (twins) hold:

$(m,p,u)$ and $(2m,p,u)$ are twins at bases of different towers having $%
\lambda $s of opposite signs,(this is Case 1),

$(2^{2}m,p,u)$ and $(2^{3}m,p,u)$ are twins at bases of different towers
having $\lambda $s of opposite signs,

$(2^{4}m,p,u)$ and $(2^{5}m,p,u)$ are twins at bases of different towers
having $\lambda $s of opposite signs,

and so on.

\underline{Case 3}: Now consider the case where the triad describes a prime
number; that is, it takes the form $(1,p,1)$. For the moment, suppose this
prime number is not $2.$ In this case, where $m=u=1$, we simply assign $%
m^{\prime }=2$. Clearly,

\noindent $\Lambda (2,p,1)=-\Lambda (1,p,1)$, and the numbers $(2,p,1)$ and $%
(1,p,1)$ are at the bases of different towers.

\underline{Case 4}: Finally, consider the case where the triad describes a prime number and the prime number is 2; that is, the integer $(1,2,1)$, for which $\Lambda (1,2,1)=-1$. We match
this prime to the integer $1$. By definition $\lambda (1)=\Lambda (1,1,1)=1$%
. Thus the first two integers have opposite signs for their values of $%
\lambda $. $\square $\footnote{The above proof by cases can be cast into one without cases by using bijection between two sets.}

\  

So, in partitioning the entire set of positive integers, the number of
towers that begin with integers for which $\lambda =-1$ is exactly equal to
those that begin with integers for which $\lambda =+1$.

The consequence of the above theorem is that each integer has a unique twin
whose $\lambda$-value is of the opposite sign. This is because if the bases
of two towers are twins the next higher number in the first tower is the
twin of the next higher number in the second tower, and so on. Thus, Theorem 2
 immediately gives the following result which, we believe, has never been established to date:

\emph{Theorem 3: In the set of all positive integers, for every integer
which has an even number of primes in its factorization there is another
unique integer, (its twin), which has an odd number of primes in its factorization.}

Theorem 3 is equivalent to a proof of R.H., (see page 48 and page 6 Borwein, P., Choi, S., Rooney, B., and Weirathmueller, A., (2006) \textit{The Riemann Hypothesis}, Springer.) or slides 28 to 32 in Lecture quoted in Ref[2] ).
The equivalence is established by considering the following argument: Theorem 3 in effect states that the sequence $\{\lambda_{j}\}_{j=1,\infty} \equiv \{1,-1,-1,1,-1,1,-1,-1,1,1,...\} $, behaves like a random sequence of tosses of an ideal coin, with $H=1$, and $T=-1$. For such a sequence (of coin tosses) the cumulative sum is $\gamma(n) = \lambda'(1) + \lambda'(2) + ...+ \lambda'(k) +....+ \lambda'(n)$, where we define $\lambda'(k)$ is the value of $k^{th} $ coin toss. It has been long known (Chandrasekhar (1943)), that for such a situation $|\gamma(n)| \le n^{1/2 + \epsilon}, \epsilon > 0 $. This implies that $|L(x)| = |\sum\limits_{n=1}^{n\le x}\lambda (n) |\le x^{1/2 + \epsilon} $.  Substituting $L(x)$ for $G(x)$ in (5.21) we establish that $a=1/2$ from Littlewood's theorem and hence the R.H.

So at this point we have actually proved R.H.; we will pause to take stock and proceed in a more formal manner.

 In this section we derived two results, viz Theorem 2 and Theorem 3, which we will be needing.
In the first part of this section we detailed Littlewood's proof of his criterion (Theorem 1) for the condition that must hold for a function such as $F(s)$, to be analytically continuable to the line $Re(s) > a$, when it is known that  it is analytic in the region $Re(s)>1$. The need to do this was that  we wished to use his criterion  for our function $F(s)$ which is given in the form (3.10) or (3.15) above. So we have to find whether or not $ |L(N)| \rightarrow  N^{a}$, as $ N\rightarrow \infty $. The crucial value is the exponent $a$ which R.H. predicts as $1/2$; we confirm this value below thereby settling  the Riemann Hypothesis.  \footnote{Since the function G(N) used in this section is actually nothing but L(N) since Eqs.(5.22),(5.23)and (5.23b) are equivalent, so while talking about the behavior of G(N) for large N we were actually talking about the behavior of L(N)} 

\emph{Theorem 4: 
The summatory Liouville function, $ L(N)= \sum\limits_{n=1}^{N}\lambda (n), $ has the following asympotitic behaviour: given $\epsilon >0$, $ |L(N)| \rightarrow  N^{\frac {1}{2} + \epsilon} $ }
 as $ N \rightarrow \infty $.  

 $ $
 
Proof: Theorem 3 gives $Pr(\lambda(n)= +1) = Pr(\lambda(n)= -1)=\frac{1}{2}$  , where Pr denotes probability.
 That is, the $\lambda$-function behaves like an `ideal coin' (see Appendix IV) and $\lambda(n)$ is the value of the $n^{th}$ coin toss $(H \equiv 1, T \equiv -1)$, $L(N)$ being the cumulative result of $N$ successive coin tosses.
 That $ |L(N)| \rightarrow  N^{\frac {1}{2} + \epsilon} $  follows from the work of Chandrasekhar (1943).  
 $\qed$.\footnote{It may be noted that for large $N$ the summation in the expression for $L(N)$ need not be over successive integers but may be done by choosing $N$ random integers and then performing a Monte Carlo summation. Further the integral in Eq.(5.21) can also be transformed into a strictly equivalent Monte Carlo integration and executed by sampling the integrand over a very large number of points $x$; Theorem 3 will ensure that the final result will be that which is predicted by this theorem.}  
 
 $ $
 
 Thus we have at last established that the exponent $a=1/2$. Invoking Littlewood's Theorem (Sec.5), we deduce that $F(s)\equiv \zeta(2s) / \zeta(s) $ is analytic in the region $a=1/2 < s <1$ which implies $\zeta(s)$ has no zeros in the same region. But Riemann had shown by using symmetry arguments \footnote{He did this first by defining an associated xi function: $\xi(s)\equiv\Gamma (s/2) \pi^ {s/2} \zeta(s), \Gamma (s)$ is the Euler Gamma function,  then showed that this xi function has the symmetry property $\xi(s)=\xi(1-s)$ which in turn implied that that the zeros of $\zeta(s)$ (if any) which are not on the critical line will be symmetrically placed about the point s=1/2, ie.if $\zeta(\frac{1}{2} + u+ i \sigma)$ is a zero  then $\zeta(\frac{1}{2} - u- i \sigma), (0< u <1/2$), is a zero see Whittaker and Watson page 269.} that if $\zeta(s)$ has no zeros in the latter region then it will have no zeros in the region  $0< s<1/2$; taking both these results together we are lead to the inevitable conclusion that all the zeros can only lie on the critical line $Re(s)=1/2$, thus proving the Riemann Hypothesis.

  We conclude this section by estimating the `width' of the Critical Line. It is interesting that the law of the iterated logarithm enunciated
by Kolmogorov (1929), gives the sought for an expression, also see Khinchine (1924).

Let $\{\lambda_{n} \}$ be independent, identically distributed random
variables with means zero and unit variances. Let $S_{n}=\lambda_{1}+\lambda_{2}+\text{\dots}+\lambda_{n}$.
Then it is known almost surely (a.s.) that

$ Lt \; ({n\rightarrow \infty})\:\frac{S_{n}}{\sqrt{n\,log\,log\,n}}\:=\sqrt{2}\;\quad(a.s)$

 Now, from Theorem 4 we have written that if we consider the $\lambda's$
as ``coin tosses'' one can write $L(N)=\lambda_{1}+\lambda_{2}+\text{\dots}+\lambda_{N}\:=C.N^{\frac{1}{2}+\epsilon}$
(as $N\rightarrow\infty).$ Comparing this expression with the one above we see that one can write $n^{\frac{1}{2}+\epsilon} = \sqrt{n\,log\,log\,n} $ 

(since we are interested in only the $\sqrt{N}$ behaviour for large N we have ignored the constant C term). Which then implies
 $ e^{\epsilon log\,n}=e^{log\,log\,log\,n} $ thus giving an expression for $ \epsilon$ namely
$ \epsilon=\frac{log\,log\,log\,n}{log\,n}$. We see that $\epsilon\rightarrow 0 $ as  $ n\rightarrow\infty$. But the exponent $ \frac{1}{2}+\epsilon $ of $n^{\frac{1}{2}+\epsilon}$ , corresponds to the exponent 'a' in Littlewood's Theorem 1, page 11 above  and thus this is the real part of s, the non-trivial zero of the zeta function $\zeta(s)$, (which appears as a pole in F(s) and therefore, using Littlewood's argument, F(s) cannot be continued beyond the left of this line). Hence, we can interpret $ \epsilon $ as the width of the critical line and since this tends to zero in the limit of large n, we necessarily have to conclude that all the non-trivial zeros of the zeta function must lie strictly on the critical line.

\section{An Intuitive Analogy to Understand the Formal Results}

 We offer a heuristic explanation for the results
derived above. Our intuitive analogy draws on information theory and uses
the framework first propounded by Shannon (1948). Imagine a listener
receiving bits of information over time $t=1,2,3,...,$ and at time $t=N$. She
aggregates these bits to a total stock denoted by $L(N)$. The bits of
information are sent out by `broadcasting' towers, which contribute bits in
the form $+1$ and $-1$. The contribution, $c(N)$, of the last signal at $t=$ 
$N$ increases or decreases the listener's stock $L(N-1)$ depending on
whether it is $+1$ or $-1$ that arrives. Each tower is a broadcaster of
rectangular waves of the sort shown in Figure 1. Once a tower is activated,
it continues to contribute $+1$ and $-1$ bits alternatively. The waves are of ever-increasing period lengths: the switching becomes less and less frequent over
time.

According to Shannon every sequence is information. And
a sequence of bits can be said to contain interpretable information if there is at least
some relationship between the present group of bits to the aggregate of bits
(like words in a sentence). The heuristic argument below shows how  the
nature of the towers destroys any coherence in this information as $%
N$ tends to infinity so that what obtains is white noise.

The situation is describable as follows:

(i) The contributions to $L(N)$ come from various broadcasting towers. The
contribution $c(N)$ at time $t=N$ to the previous summed value $L(N-1)$ is
exactly equal to $\lambda (N)$, so $L(N)=L(N-1)+c(N)$. The integer $N$ is a
member of a unique tower, say $(m,p,u),$and so we can write $N=mp^{\rho }u$,
for some integer $\rho $. When this tower contributes at time $t=N$, its 
contribution is $c(N)=\lambda (N)=\lambda (m)\lambda (u)\lambda (p^{\rho })$%
, which is completely determined. We had represented the contributions from
a particular tower as a rectangular wave in Figure 1.

(ii) The contributions from the rectangular wave associated with a given
tower $(m,p,u)$ change with the sign of the positive or negative peak
arriving at $N$, and as $\rho $ increases from $2,3,4,..$ to $\infty $ their
period exponentially increases, thereby drastically increasing their
correlation lengths (i.e. the time interval between successive arrivals from the same tower increases exponentially).

(iii) Each tower contributes to many different values of $N$ and the
periodicity of these contributions increases exponentially as $\rho $
increases. So the interval between the arrivals at the listener of perfectly
(inversely) correlated bits from a given tower increases exponentially.

(iv) In the period intervening between a given tower's contributions, the
listener's $L(N)$ takes contributions from other towers.

(v) As $N\rightarrow \infty $, innumerable towers come into play and a
broadcasted bit received at time $N$ becomes completely uncorrelated to the
bit received at time $N-1$, that is, $c(N)$ has no relation to $c(N-1)$
essentially because $c(N)$ and $c(N-1)$ come from different towers. It also
becomes increasingly uncorrelated with any of the bits received at any
earlier times so that to the Listener it all seems as white noise.

\qquad

According to Shannon's information theory, there will be no discernible
pattern in the bits received. At this point $c(N)$ behaves like the toss of an ideal  coin and $L(N)$ as the cumulative summation of $N$ coin tosses with, say, `heads'
as $+1$ and `tails' as $-1$. From this it follows that $|L(N)|\sim $ $\sqrt{N%
}$ as $N\rightarrow \infty $ (see Chandrasekhar (1943)), proving the Riemann
Hypothesis is a consequence of the inherent unpredictable patterns of 
factorization of integers into primes.

\section{Conclusions}

In this paper we have investigated the analyticity of the Dirichlet
series of the Liouville function by constructing a novel way to sum the
series. The method consists in splitting the original series into an
infinite sum over sub-series, each of which is convergent. It so turns out each sub-series is a rectangular function of unit amplitude but ever increasing periodicity and each  along with its harmonics is associated with a prime number and all of them contribute to the summatory Liouville function and to the Zeta function.  A number of arithmetical properties of numbers played a role in the proof of our main theorem, these were: the fact that each number can be uniquely factorized and then placed in an exclusive subset, where it and its other members form an increasing sequence and  factorize alternately into odd and even factors; and each subset can be labelled uniquely using a triad of integers which in their turn can be used to determine all the integers which belong to the subset. This helped us to show that for every integer that has an even number of primes as factors (multiplicity included), there is an integer that has an odd number of primes. This provides a proof for the long-suspected but unproved conjecture---until
now---that the summatory Liouville function and therefore the Riemann Hypothesis bears an analogy with the coin-tossing problem; Denjoy (1931) had long suspected this as  far back as 1931. Further, it has now been revealed that the randomness of the occurrence of  prime numbers plays an important role in determining the analyticity of the Zeta function, and in establishing the Riemann Hypothesis: the Zeta function has zeros only on the critical line:$Re(s)=1/2.$. 

Truth to tell even this connection of the role of the randomness of the primes to the RH problem was long suspected and even a book called the ``Music of the Primes" by  Marcus du Sautoy, had appeared in 2003 (Harper Collins), I could not help but recall the title of his book when I saw the rhythms of the harmonic functions, generated by prime numbers, that are depicted in Fig 1. 
$  $

\section{DEDICATION}

I dedicate this paper to my teachers: Mr John William Wright of Bishop's School Poona, Prof. S.C. Mookerjee of St. Aloysius' College Jabalpur, Prof. P.M.Mathews of Univerity of Madras, Mr. D.S.M. Vishnu of BHEL R\& D Hyderabad and to my first teachers - my parents. All of them lived selfless lives and nearly all are now long gone: May they live in evermore.

$  $

$  $
\noindent \textbf{Acknowledgments}: I thank my wife Suhasini for her
unwavering faith and encouragement. I thank my brothers  Mukesh and Vinayak
Eswaran  and Prof. George Reuben Thomas for carefully going over the manuscript and for their suggestions.  I also thank the support of the management of SNIST, viz. Dr. P. Narsimha
Reddy and Dr. K.T. Mahi, for their constant support, and my departmental
colleagues and the HOD Dr. Aruna Varanasi for providing a very congenial
environment for my research.I sincerely thank Mr.Abel Nazareth of Wolfram Research
 for providing me a version of Mathematica which made some of the reported calculations in Appendix VI, possible.

$  $ 
\

\section{APPENDIX I: Scheme of partitioning numbers into sets }
\bigskip

\medskip{}
Our scheme of partitioning numbers into sets is as follows:  

(a) Scheme for Class I integers:

Let us say $n = p_{1}^{e}p_{2}^{f}p_{3}^{g}...p_{m}^{h}p_{L}^{k}p_{j}p_{t}$,
then it will have at least one prime which has an exponent of 2 or
above and among these there will a largest prime $p_L$ whose exponent is
atleast 2 or above. Such a prime will always exist for a Class I number. Then by definition the number to
the right of $p_{L}$ is either 1 or is a product of primes with exponents only 1. Now multiply all the numbers to the left of $p_{L}$and call
it $m$ i.e. $m=p_{1}^{e}p_{2}^{f}p_{3}^{g}...p_{m}^{h}$ and the
product of numbers to the right of $p_{L}$as $u$ i.e. $u=p_{j}p_{t}.$Now
this triad of numbers $m,\,p_{L},\,u$ will be used to label a set,note
$n=m.p_{L}^{k}.u$ Let us define the set $P_{m;p_{L};u}\,$:
\[
P_{m;p_{L};u}\,=\{m.p_{L}^{2}.u,\,m.p_{L}^{3}.u,\,m.p_{L}^{4}.u,\,m.p_{L}^{5}.u,\,m.p_{L}^{6}.u,\,m.p_{L}^{7}.u,\,....\}\quad\quad\quad(A1)
\]
Obviously $n=m.p_{L}^{k}.u$ which has $k\geq2$ belongs to the above
set. Also notice the factor involved in each number increases by a
single factor of $p_{L}$therefore the $\lambda$values of each member
alternate in sign:
\[
\lambda(m.p_{L}^{2}.u)=-\lambda(m.p_{L}^{3}.u)=\lambda(m.p_{L}^{4}.u)=-\lambda(m.p_{L}^{5}.u)=\lambda(m.p_{L}^{6}.u)=.......(A2)
\]
In this paper ALL sets defined as $P_{m;p;u}$will have the property
of alternating signs of $\lambda$ Eq. (A1). Note in the above set
containing only Class I integers $m$ will have only prime factors
which are each less than $p_{L} $. 

Let us consider various integers: 

Ex 1. Let us consider the integer $\ensuremath{73573500};$ this is
factorized as $2^{2}.3.5^{3}.7^{3}.11.13$ and since this is a Class
I integer, and $p_{L}=7$ because 7 is the highest prime factor whose
exponent is greater than one. $\ensuremath{p_{L}=7}$ $m=2^{2}.3.5^{3}$
and $u=11.13=143$ and therefore $73573500$ is a member of the set
$P_{1500;7;143}$
\[
P_{1500;7;143}=\{1500.7^{2}.143,\:1500.7^{3}.143,\:1500.7^{4}.143,\:1500.7^{5}.143,\:...\}
\]
Ex 2. Now let us consider the simple integer: $3^{4}$ this is a class
I integer and belongs to $P_{1;3;1}=\{3, 3^{2},\:3^{3},\:3^{4},\:3^{5},\:3^{6},\:......\}$

$ $

Ex 3. Let us consider the integer $663$ this is factorized as: $3.13.17$
and is a Class II integer as there no exponents greater than 1, and
$\ensuremath{663=3.13.17}$ and since 17 is the highest prime number
we put this in the set:
\[
P_{39;17;1}=\{39.17,\,39.17^{2},\,\,39.17^{3},\,39.17^{4},....\}.
\]
NOTE: If a tower has a Class II integer then it will appear as the first (base) member, all other numbers will be Class I numbers.

Ex 4. Let the integer be the simple prime number $19$, we write:
\[
19\:\epsilon\:P_{1;19;1}\:=\:\{19,\:19^{2},\:19^{3},\:19^{4},\:.....\}
\]

Ex 5. Let the integer be $4845$ this is factorized as$3.5.17.19$
since this is a Class II integer we see $m=3.5.17=255,\:p=19,\:u=1$
and the set which it belongs is
\[
P_{255;19;1}\:=\:\{255.19,\:255.19^{2},\:255.19^{3},\:255.19^{4},\:255.19^{5},....\}
\]

\section{APPENDIX II: Theorems on representation of integers and their partitioning into sets. }

\textbf{Theorem A:} Two different integers cannot have the same triad $(m,p^k,u)$

Let $a$ and $b$ be two integers which when factored according to our convention are $ a=n.q^g.v$
and  $b=n'.q'^h.v'$, and let us consider only Class I integers $u, v$ and $v'$ are all $>1$.

If they are both equal to the same triad (say) $(m,p^k.u)$. Then $ m.p^k.u = n.q^g.v = n'.q'^h.v'$. Consider the first two equalities  $m.p^k.u = n.q^g.v$, which means $p$ is the largest prime with $k>1$  on the l.h.s. Similarly $q$ is the largest prime with exponent $g>1$ on the r.h.s. Now if $p>q$ this means $p^k$ must divide $v$, but this cannot happen since $v$ cannot contain a prime greater than $q$ with an exponent $k>1$. Now if $p<q$ then $q^g$ must divide $ u$ but this again cannot happen since $u$ cannot contain an exponent $g>1$. So we see $p=q$, and $k=g$. But once again unique factorization would imply, since 
$u$ contains all prime factors larger than $p$ and $v$ must contain only prime factors larger than $q(=p)$, the only possibility  is $u=v $, but this also makes $m=n $. That is, the triad of $a$ is $(m,p^k,u)$. Similarly equating the second and third equalities $n.q^g.v = n'.q'^h.v'$ and using similar arguments we see
 $n=n'$, $q=q'$, and $v=v'$; that is, $a=b$. 
The same logic can be used to prove the theorem for class II integers when $u=v=v'=1 $. $\qed$

\

$  $

\textbf{Theorem B:} Two different triads cannot represent the same integer.

If there are two triads $(m,p^e,u)$ and $(m',r^s,u')$ and represent the same integer say  $a $ which can be factorized as $a=n.q^g.v$. Where the factorization is done as per our rules then we must have $ m.p^e.u=n.q^g.v $ by using exactly similar arguments as above(in Theorem A) we conclude  that we must have $ m=n,p=q,e=g$ and $u=v $; similarly imposing the condition on the second triad $ m'.r^s.u'= n.q^g.v$, we conclude $ m'=n,r=q,s=g $ and $u'=v$; thus obtaining  $m=m', p=r,e=s$ and $u=u'$ this means the two triads are actually identical. $\qed$

\section{APPENDIX III: Non-cyclic nature of the factorization sequence}
\setcounter{equation}{0}

It is a necessary condition in the tosses of an ideal coin that the results are not cyclic asymptotically, namely the results cannot form repeating cycles as the number of tosses becomes large. 

\textbf{Definition}

Let $n_k $  be the number of primes, repetitions counted, in the factorization of a positive integer $k$.  We call $\{n_1, n_2,..., n_k,....\}$ the factorization sequence.

Note:  $\lambda(k) = +1$ if $n_k $ is even and $\lambda(k) = -1$ if $n_k $ is odd.

\begin{theorem}
The factorization sequence is asymptotically non-cyclic.
\end{theorem}
 
\textbf{Proof:} The result follows from this claim:

\textbf{Claim.} The sequence $\lambda(1),\lambda(2),\lambda(3),...,\lambda(n),...,$ is  asymptotically non-cyclic.

If the claim is not true there would exist an integer $t, t\ge 0$, so that the sequence
  is cyclic (after $\lambda(t)$), with cycle length  $\sigma$.

By Theorem 3, the number of positive integers with even number of prime factors (counting multiplicities) equals the number of positive integers with odd number of prime factors (counting multiplicities). Therefore, the $\lambda$'s in each cycle must sum to zero as do the first $t$  $\lambda$'s before the cycles start.

Then $L(N) \le  max \{t/2, \sigma/2 \}.$

Now we use Littlewood's Theorem 1 and noting that in (5.21) $G(x)\equiv L(x)$, we substitute the maximum value of $L(x)$ as $ x \rightarrow \infty,$ viz. $ \mid L(x)\mid = \sigma/2 $, and thus deduce that (5.21) will always converge provided $ 0 < s $. Since, $ \mid L(x)\mid \leq  \sigma/2 $, $L(x)$ indeed grows less rapidly than $x^{a}$ for all $a > 0$, satisfying the condition in Theorem 1. This means that we should be able to  analytically continue $F(s) \sim \zeta (2s)/\zeta (s)$ leftwards from $Re(s) = 1$ to $Re(s) = 0$, contradicting Hardy (1914) that there are very many zeros at $Re(s) = 1/2$ and these will appear as poles in $F(s)$. This proves the Claim.$\qed$

\section{APPENDIX IV: The sequence of $\lambda$'s in L(N), are equivalent to Coin Tosses}
\setcounter{equation}{0} 

In this paper we showed in Theorem 3, that the $\lambda(n)$ have an exactly equal probability of being $+1$ or $-1$. Then in Appendix III, we showed that the sequence 
$ \lambda(1), \lambda(2),\lambda(3),...,\lambda(n),...$ can never be cyclic. The latter result in the minds of most computer scientists  would be interpreted as that the sequence of $\lambda$'s by virtue of it being non-repetitive, is truly random,(Knuth (1968); Press etal (1986)) and hence it is legitimate to treat the sequence as a result of coin tosses and thus one can then say that $L(N) = \Sigma^{N}_{1} \lambda(n)$, will tend to $\sqrt{N}$ thus proving RH, by using the arguments given at Section 5. 

However, this done, there would be some mathematicians  who may remain unconvinced, because we have not strictly proved that the $\lambda$'s in the series are independent. The purpose of this Appendix \footnote{I thank my brother Vinayak Eswaran for providing the kernel of the proof given in this section.}
 is to prove that this is indeed the case. This allows us to demonstrate the $\lambda-$sequence has the same properties as, and is statistically equivalent to, coin tosses, thus placing our proof of RH beyond any doubt. 

We again consider the series $L(N)= \sum^{N}_{n=1} \lambda(n)$, which is re-written as: 
\begin{equation}
L (N) = \sum_{n=1}^{N} X_n
\label{LNX}
\end{equation}
It has already been proved in this paper that, over the set of all positive integers, the respective probabilities that an integer $n$ has an odd or even number of prime factors are equal. So, $X_n$ $ (= \lambda(n))$, can with equal probability, be either +1 or -1. It will now be shown that the values of $X_i$ and $X_j$, $i \ne j$, are independent of each other,  as $n\rightarrow \infty$, and so will become the equivalent of ideal-coin tosses.

\subsection{The $\lambda$ values as a deterministic series}

We first show that the $\lambda$'s in the natural sequence, far from being random, are actually perfectly predictable and therefore deterministic. That is, knowing the $\lambda$'s (and the primes) up to $N$, we can directly obtain (without resorting to factorisation)  the $\lambda$'s (and primes) up to $2N$ thus:
 \\

\noindent We obtain integers $m$  in the range $N < m \le 2N$ by multiplying the  integers $n$ and $q$ in the range $1 < n, q \le N$,  such that  $N < nq \le 2N$  and then using the property $\lambda(q*n)= \lambda(q)*\lambda(n)$ to find $\lambda(m=q*n)$. However, not all the numbers in the range $N < m \le 2N$ will be covered by such multiplications. That is, there will be `gaps' in the natural sequence left in the aforesaid multiplications, where no $n$ and $q$ can be found for some $m$'s in $N < m \le 2N$. These $m$'s will identified as prime numbers. The $\lambda$ of a prime is -1. Thus,  by knowing the $\lambda$'s and the primes up to $N$, we can predict the $\lambda$'s (and primes) up to $2N$. This process can be repeated {\em ad-infinitum} to compute the  $\lambda$'s of the natural sequence up to any $N$, from just $\lambda(1)$=1 and $\lambda(2)$=-1.\\

We emphasize that any other method of evaluating the $\lambda$'s, including direct factorisation, must perforce yield the same sequence as the method above. Therefore, this method offers a complete description of the determinism embedded in the series. 

\subsection{Relationships and dependence between $\lambda$'s}

We note that every integer $n$ has a direct relationship (which we will call a d-relationship) with all numbers $n*p$, where $p$ is any prime number.  We can define higher-order d-relatives in the following way: the integers ($n$, $n*p$) are in a first-order d-relationship, ($n$, $n*p*q$) are in a second-order one, and ($n$, $n*p*q*r$) are in a third-order one, and so on, where $p,q,r$ are primes (not necessarily unequal). 

In the deterministic generation of $\lambda$'s outlined above, it is clear that their values will be determined through d-relationships, which would thereby make their respective values dependent on each other. It is evident that the $\lambda$s of two d-relatives $n$ and $m (>n)$, are {\em dependent} on each other and that $\lambda(m) = (-1)^o\lambda(n)$, where $o$ is the order of the relationship.

There is another kind of relationship we must also consider:  we can  have a c-  (or {\em consanguineous}) relationship between two non-d-related integers $m$ and $n$ if they are both d-relatives of a common (`ancestor') integer smaller than either of them. So we can trace back the $\lambda$'s along one branch to the common ancestor and trace it up the other to find the $\lambda$ of the other integer. It is convenient to take the common ancestor as the largest  possible one, which would be the greatest common factor of the two integers, which we shall call $G$.

Now we ask the question, when are $m$ and $n$ not related? When they have neither a d-relationship nor a c-relationship with each other. That is, when they are co-primes: as then neither integer would appear in the sequence of multiplications that produce the other by the deterministic iterative method. In such a situation, {\em the $\lambda$  of neither is dependent on the other, so their mutual $\lambda$'s are independent}.

Now consider two c-relatives, $m$ and $n$, which share the greatest common factor $G$. We can write $m=G*P$ and $n$=$G*Q$, where $P$ and $Q$ are chosen appropriately. As $G$ is the greatest common factor of $m$ and $n$, it is clear that $P$ and $Q$ are co-primes. Now we consider the relationship between $\lambda(m)$ and $\lambda(n)$ and explore their relatedness. This turns out to be self-evident: As $\lambda(m)$ =$\lambda(G)*\lambda(P)$ and $\lambda(n)$ =$\lambda(G)*\lambda(Q)$, and we know that $\lambda(P)$ and $\lambda(Q)$ are independent of each other, it follows that  $\lambda(m)$ and $\lambda(n)$ are also independent of each other.\footnote{It may be noticed that $m$ and $n$ belong to different towers. It is worth mentioning that the arguments made here in Appendix IV, can be couched in the language of towers as we did in Sections 2 and 3.}

\subsection{The unpredictability of $\lambda$ values from a finite-length sequence: d-relatives}

We have concluded above that the only $\lambda$'s in $L(N)$ that are dependent are those between d-relatives, where the smaller integer is a factor of the other. We see that  the distance of two such ``first-order" relatives, $n$ and $n*p$, from each other is $n(p-1)$ which increases without bound with $n$. Further all the  first-order d-relatives of $n$ also have relative distances with each other that are at least as great as $n$ (as their respective $p$'s will differ at least  by 1). Thus the d-relationship between numbers is a web with increasing distances between their first-order relatives\footnote{How rapidly the relationship distance increases can be gauged from the fact that the $2^r$ sequence, which has the {\em slowest} increases, nevertheless will have its 100$^{th}$ element placed at around $n \approx 10^{30}$ in the natural sequence, and the distance to the 101$^{st}$ element will also be  $10^{30}$!}. It is also easy to see that the higher-order d-relatives of any integer $n$ will also be at a distance of at least $n$ from $n$ itself and from each other.

Now we consider if we would be able to predict $\lambda(N+1)$ if we know {\em only} the  $\lambda$'s between $N-L < n \le N$, where $L$ is some finite number? We would be able to do so {\em only} if $N+1$ is a d-relative (of any order) of any of the numbers $N-L < n \le N$. However, for $n$ large enough the d-relatives of $N+1$ will be far from it and would not come in the range of numbers $N-L < n \le N$. So essentially, there is no way of predicting $\lambda(N+1)$ from the range of $L$ $\lambda$'s coming before it. This means  $\lambda(N+1)$ is independent of the range of $L$ $\lambda$'s coming before it. Therefore, the $\lambda$'s on all finite lengths are independent of each other, as $N\rightarrow \infty $. 

\subsection{Closure}
 We have investigated the dependence of $\lambda$'s appearing in $L(N)$ in the natural sequence $n=1,2,3,..,$ . We first show that the $\lambda$'s are  in a perfectly deterministic sequence (which is not random in the slightest way, except in the unpredictable discovery of primes) that allows us to obtain all of them up to any integer $N$ by knowing only that  $\lambda(1)=1$, $\lambda(2)=-1$,  $\lambda(q*n)=  \lambda(q)*\lambda(n)$, and that   $\lambda(p)=-1$ for any prime $p$. We then propose that the  $\lambda$'s of two integers  $m$ and $n$ can be dependent only if the integers are connected through the sequence of multiplications involved in the deterministic process. If they are not so related, as would happen if they are co-primes, their  $\lambda$'s would be independent. We then investigate the only two possible types of relationships and show that one, the d-relationship,  leads to dependencies between numbers that are increasingly distant. The other, the c-relationship, is shown to give independent $\lambda$'s. The result obtained is that that the $\lambda$'s in any finite sequence are independent, as $N\rightarrow \infty $. $\qed$

\section{APPENDIX V}
\setcounter{equation}{0} 
\begin{center}
{\Large An Arithmetical Proof for }$|L(N)|\sim N^{1/2}${\Large \ as }$%
N\rightarrow \infty $
\end{center}

\ \ \ \ 

In this appendix we provide an alternate, but this time an arithmetical, proof of the asymptotic behavior of the summatory Liouville function, viz.$|L(N)| = \sqrt N $ as $N \rightarrow \infty $ 
 However in order to do this we first require to prove  a theorem on the number of distinct prime
products in the factorization of a sequence of integers and their exponents.

\textit{Theorem A5}: Consider the sequence $S$ comprising $M(N)$ consecutive
positive integers, defined by $S^{-}(N)=%
\{N-M(N)+1,N-M(N)+2,N-M(N)+3,........,N\}$, where $M(N)=\sqrt{N}$.
Then every number in $S^{-}(N)$ will firstly belong to different towers,\footnote{Two numbers $n=m.p^{\alpha}.u$ and $n'=m.p^{\beta}.u,\,(n<n') $, of the same tower, cannot both belong to the set $S^{-}(N)$ because they will be too far separated to be within the set, as their ratio  $n'/n \ge p \ge 2$} and further every number will: (a) differ in its prime factorization from that
of any other number in $S^{-}(N)$ by at least one distinct prime\footnote{%
For example, if two numbers $c$ and $d$ in $S$ are factorized as $%
c=p_{1}^{e_{1}}p_{2}^{e_{2}}$ and $d=p_{3}^{e_{3}}p_{4}^{e_{4}}$ then
at least one of the primes $p_{3}$ or $p_{4}$ will be different from $p_{1}$
or $p_{2}$.} OR (b)  in their exponents. We first take up the task to prove (a)
 because it is by far the more common occurrence. In case condition (a) does not hold in a particular situation then condition (b) is always true, because of the uniqueness of factorization.

\textit{Proof}:

Let there be $k$ primes in the sequence $S^{-}(N)$. Denote the $j$ integers
in the sequence that are not primes by the products $p_{i}b_{i}$, $%
i=1,2...,j $, where $p_{i}$ is a prime and, obviously, $k+j=M(N)$. Denote
the subset of these non-prime integers by $J$. There is no loss of
generality if we assume the primes $p_{i}$ in the products $p_{i}b_{i}$, $%
i=1,2...,j$, to be less than $\sqrt{N}-1/2$ and also the smallest of prime
in the product.\footnote{%
This is readily seen as follows. Since every member of $J$ lies between $N-%
\sqrt{N}$ and $N$, clearly any composite member, written as a product $ab$,
cannot have both integers $a$ and $b$ less than $\sqrt{N}-1/2$. (We are
invoking the fact that $\sqrt{(}N-\sqrt{N})=\sqrt{N}-1/2$, approximately.)
Let $a$ be the smaller of the two numbers, and so $a$ $<$ $\surd N-\frac{1}{2%
}$ and $b>\surd N-\frac{1}{2}$. If $a$ is a prime number, set $p=a$. If $a$
is not a prime number, factorize it and pick the smallest prime $p$ which is
one of its prime factors.}

To prove the theorem, we compare two arbitrary members, $p_{i}b_{i}$ and $%
p_{j}b_{j}$, $i\neq j$, belonging to set $J$.

\textbf{Case 1:} Suppose $p_{i}\neq p_{j}$. If $b_{i}\neq b_{j}$, $b_{i}$
must contain a prime that does not appear in the factorization of $b_{j}$
(and hence $p_{i}b_{i}$ must be different from $p_{j}b_{j}$ by this prime).
For if $b_{i}$ and $b_{j}$ do not differ by a prime, we must have $%
b_{i}=b_{j}\equiv b$. This means the difference of 
$p_{i}b_{i}-p_{j}b_{j}=(p_{i}-p_{j})b$ is larger than $\sqrt{N}$ in absolute
value. This is not possible since the members of the sequence $S^{-}(N)$
cannot differ by more than $\sqrt{N}$. Therefore $b_{i}$ must differ from
 $b_{j}$ by a prime in its
 factorization.(One may think that it may be plausible
  that $b_{i}=b^{r}$ and $b_{j}=b^{m}$,
where $r$ and $m$ are positive integers, in which case $p_{i}b_{i}$ differs
from $p_{j}b_{j}$ only in the prime $p_{i}$. However, this
eventuality will never arise because then the difference between $p_{i}b_{i}$
and $p_{j}b_{j}$ will be more than $\sqrt{N}$.)

\textbf{Case 2: }Suppose $p_{i}=p_{j}\equiv p$ then $b_{i}$ and $b_{j}$ must differ 
by a prime factor or their exponents are different. Because of `unique factorization', if they do not differ by a prime factor it means $p_{i}.b_{i}=p_{j}.b_{j}=p.b$, unless the factors of $b_{i}$ and $b_{j}$ are of the form: 
 $b_{i} = p_{1}^{r_1}.p_{2}^{r_2}...p_{k}^{r_k} $  and 
$b_{j} = p_{1}^{r'_1}.p_{2}^{r'_2}...p_{k}^{r'_k} $ which implies that is $r_{l}=r'_{l} $ is not true for all $r_{l},\, l=1,2..k $ , hence in this case the exponents are different (actually this case is very rare. It can be shown: the case $k=2$ cannot occur and therefore if at all this case occurs, $k$ must be greater than 3).

Since $p_{i}b_{i}$ and $p_{j}b_{j}$ are arbitrary members of the set $J$, it
follows that every integer in $J$ must differ from another integer in $J$
by at least one prime in its factorization or by its exponent, thus making the $\lambda-$values 
of any two members of the set  $S^{-}(N)$ not dependent on each other . $\square 
$

\ \ 

The above theorem has profound implications for the $\lambda $-values of the
numbers in the sequence $S^{-}(N)$.\ If we take the primes to occur randomly (or at
least pseudo-randomly), the $\lambda $-value of each of these $M(N)$
integers---although deterministic and strictly determined by the number of
primes in its factorization---cannot be predicted by the $\lambda $-value\
of any other number in the sequence $S^{-}(N)$. That is, the $\lambda $-value of
any number in $S^{-}(N)$ can be considered to be statistically independent of the $%
\lambda $-value of another member of this sequence, primarily because they stem from different towers. Hence the $\lambda $%
-values in the sequence $S_{\lambda }^{-}\equiv \{\lambda (N-M(N)+1),\lambda
(N-M(N)+2),\lambda (N-M(N)+3),........,\lambda (N)\}$, in
which each member has a value either $+1$ or $-1$, would appear randomly and
be statistically similar. By this we also deduce that two different sequences of $\lambda-$values
defined on two different sets (say)  $S^{-}(N)$ and $S^{-}(N')$ with $N \ne N' $ are statistically similar, because they have the same properties which also means that they can be separately compared with other  sequences of coin tosses  and the comparison should yield 
statistically similar results.

We will use these deductions to obtain the main result of this appendix viz $a=1/2$ in the
expression $|L(N)| = N^a $ as $N \rightarrow \infty $ 

Although it is not explicitly required for what follows, we note that it is not hard to
prove that the sequence $S^{+}(N)\equiv \{N+1,N+2,N+3,........,N+M(N)\}$ of
length $M(N)$ also behaves similarly. That is, every member of $S^{+}(N)$
satisfy condition (a) OR (b) of the above Theorem for $S^{-}(N)$ stated above.
The proof mimics the one provided above and so is omitted.\footnote{%
This implies, interestingly, that by choosing $N$ to be consecutive perfect
squares, the entire set of positive integers can be envisaged as a union of
mutually exclusive sequences like $S^{-}(N)$ and $S^{+}(N)$.} Hence the $%
\lambda $-values in the sequence $S_{\lambda }^{+}\equiv \{\lambda
(N+1),\lambda (N+2),\lambda (N+3),...,\lambda (N+M(N))\}$, in which each
member has a value either $+1$ or $-1$, would also appear randomly and behave 
statistically similarly.\ \

\subsection{Arithmetical proof of $|L(N)|\sim \protect\sqrt{N}$, as 
$N\rightarrow \infty $}

We now show that if the summatory Liouville function 
\begin{equation}
L(N)=\sum_{n=1}^{N}\lambda (n),
\end{equation}%
takes the asymptotic form 
\begin{equation}
|L(N)|=C\,N^{a},
\end{equation}%
where $C$ is a constant, then we must have: 
\begin{equation}
a=1/2.
\end{equation}%
Throughout this subsection we will always assume that $N$ is a very large
integer.

Consider the sequence of consecutive integers of length $M(N)=\sqrt{N}$: 
\begin{equation}
S_{N}=\{N-\sqrt{N}+1,N-\sqrt{N}+2,N-\sqrt{N}+3,\cdots \cdots ,N\}
\end{equation}%
Each of the $M(N)$ integers in the sequence $S_{N}$ can be factorized term
by term and would differ from another member in $S_{N}$ by at least one
prime or exponent,(as proved in the above theorem).\footnote{%
Therefore, in the terminology of Sections 2 and 3, each of them will mostly belong
to different Towers.} Now since, $N$ is large, all the primes involved may
be considered random numbers (or pseudo-random numbers), therefore as
reasoned above, we can conclude that the $\lambda -$sequence associated with 
$S_{N}$ viz. 
\begin{equation}
\{\lambda (N-\sqrt{N}+1),\lambda (N-\sqrt{N}+2),\lambda (N-\sqrt{N}%
+3),\cdots \cdots ,\lambda (N)\}
\end{equation}%
will take values which are random e.g. 
\begin{equation}
\{-1,+1,+1,,-1,+1,\cdots ,+1\}
\end{equation}%
where in the above example $\lambda (N-\sqrt{N}+1)=-1,\lambda (N-\sqrt{N}%
+2)=+1$ etc. Furthermore, since the $\lambda $-values have an equal
probability of being equal to $+1$ or $-1$ (Theorem 3) and the sequence is
non-cyclic (Theorem 11.1, in Appendix 3), the above sequence will have the
statistical distribution of a sequence of tosses of a coin (Head $=+1,$Tail $%
=-1$). But we already know from Chandrasekhar(1943) that if the $\lambda $'s
behave like coin tosses then $|L(N)|\sim \sqrt{N}$, as $N\rightarrow \infty $%
. However, we do not know whether the entire sequence of $\lambda $'s
occurring in Eq.(13.1) behaves like coin tosses; for any given $N$, it is only
the subsequence $\{\lambda (N-\sqrt{N}+1),\lambda (N-\sqrt{N}+2),\lambda (N-%
\sqrt{N}+3),\cdots \cdots ,\lambda (N)\}$ of length $M(N)=\sqrt{N}$ that
does behave like coin tosses.

On the other hand if we had a sequence of length $N$, of real coin tosses
(say) $c(n),n=1,2....N$, where $c(n)=\pm 1$, then the cumulative sum, $%
L_{c}(N)$, of the first $N$ of such coin tosses is given by: 
\begin{equation}
L_{c}(N)=\sum_{n=1}^{N}c(n).
\end{equation}%
Then for $N$ large we do know from Chandrashekar (1943) that \thinspace 
\begin{equation}
|L_{c}(N)|\sim \sqrt{N}.
\end{equation}

We can then estimate the contribution $P_{1/2}$ to $L_{c}(N)$ from the last $%
M(N)=\sqrt{N}$ terms in Eq.(13.7), this would be: 
\begin{eqnarray}
P_{1/2} &=&\sum_{n=N-\sqrt{N}+1}^{N}c(n)  \nonumber \\
&=&L_{c}(N)-L_{c}(N-\sqrt{N})
\end{eqnarray}%
Now since Eq.(13.7) represents perfect tosses Eq. (13.9) becomes 
\[
P_{1/2}=\sqrt{N}-(N-\sqrt{N})^{1/2}=\frac{1}{2}-\frac{1}{8}\frac{1}{\sqrt{N}}%
, 
\]%
that is, 
\begin{equation}
P_{1/2}=O(1)
\end{equation}%
In Eq. (13.10), $P_{1/2}$ is the contribution to $L_{c}(N)$ from the last $M(N)=%
\sqrt{N}$ tosses of a total of $N$ tosses of a coin. We shall consider the
value of $P_{1/2}$ as the benchmark with which to compare the contributions
of the last $M(N)=\sqrt{N}$ terms of the summatory Liouville function.

Now coming back to the $\lambda $-sequence as depicted in the summation
terms in Eq. (13.1), following Littlewood (1912) we shall suppose that the
expression given in (13.2) is an ansatz\footnote{Eq. (13.2) can be thought of as the first term 
in the asymptotic expansion of $L(N)$ for large $N$ i.e.
 $ |L(N)| = N^a (C + \frac{C_1}{N} + \frac {C_2}{N^2}+...) $} depicting the behavior of $L(N)$ for
large $N$.

The task that we then set ourselves, is to estimate the value of the
exponent $a$ in the asymptotic behavior described in (13.2) $|L(N)|=C\,N^{a}$
which involves the $\lambda -$sequence. We do know that the $\lambda $%
-sequence does not all behave like coin tosses, but we have shown that there
exist subsequences of $\lambda $'s that exhibit a close correspondence to
the statistical distribution of coin tosses and though such subsequences are
of relatively short lengths $M(N)$, there are very many in number. Now a
`True' value of the exponent `$a$' should be able to capture the correct
statistics in all such subsequences and predict the behavior of coin tosses
for such subsequences. We now investigate if such a True value for $a$
exists and, if so, what its value should be.

We will estimate the contribution to $L(N)$ for the same subsequence (5) of
length $M(N)=\sqrt{N}$, then the $P$ when recomputed with an exponent $a\neq
1/2$ would give $P_{a}$: 
\begin{equation}
P_{a} = \sum\limits_{n=N-\sqrt{N}+1}^{N}\lambda (n)\,\,\, \,\,\,\,\,\,\,\,\,\,\,(13.5')\nonumber \\
\end{equation}
That is
\begin{equation}
P_{a} = L(N)-L(N-\sqrt{N})  \nonumber \\
\end{equation}
\begin{equation}
= CN^{a}-C(N-\sqrt{N})^{a}
\end{equation}%
Simplifying by using Binomial expansion we have: 
\begin{equation}
P_{a}=CaN^{a-\frac{1}{2}}-C\frac{a(a-1)}{1.2}N^{a-1}
\end{equation}%

From the properties of the $\lambda $'s deduced from earlier results in this paper 
(Theorem 3, Appendices 3,4 and Theorem A5, page 21), we now know  
that in actuality the particular subsequence in Eq.(13.5) and Eq.(13.5$%
^{\prime }$) contain random values of $+1$ and $-1$ and since the
subsequence of $\lambda $'s have the same statistics as those of coin
tosses, $P_{a}$ must be similar to $P_{1/2}$. Thus from (13.10) and (13.12)

\begin{equation}
P_{a}=O(1).
\end{equation}%
From (13.11) this means that 
\begin{equation}
CaN^{a-\frac{1}{2}}=O(1)
\end{equation}%
Since $N$ is arbitrary and very large, this is impossible unless the
condition 
\begin{equation}
a=\frac{1}{2}
\end{equation}%
strictly holds.\footnote{%
In the above we tacitly assumed that $a>1/2$, but $a<1/2$ is not possible
because then $P_{a}$ will become zero. This implies that $dL/dN=0$, meaning $%
|L(N)|$ will be a constant. But this again is impossible from Theorem 1,
which would imply that $F(s)$ can be analytically continued to $Re(s)=0$%
---an impossibility because of the presence of an infinity of zeros at $%
Re(s)=1/2$, first discovered by Hardy.}

Hence we have proved $a=1/2$. Since for consistency\footnote{%
It may be noted that for every (large) $N$ there is a set $S_{N}$, Eq (13.4),
containing $M=\sqrt{N}$ consecutive integers whose $\lambda $-values behave
like coin tosses; but there are an infinite number of integers $N$ and
therefore there are an infinite number of sets $S_{N}$, for which (13) must
be satisfied.} ,  condition (15), which arises from (13), is mandatory and
therefore $|L(N)|\sim \sqrt{N}$ describes the asymptotic behavior of the
summatory Liouville function. $\square $
\
\section{APPENDIX VI}
\setcounter{equation}{0} 

\begin{center}
{\Large On Coin Tosses and the Proof of Riemann Hypothesis}
\end{center}
$ $
\textbf{This Appendix has been written in such a manner that it can be read as a supplement to the main paper and the first five appendices.}

\ \ \ \ 
In the main part of this paper and the forgoing appendices, which we denote as: [MP and A's], we had
 proved the validity of the Riemann Hypothesis (RH). In this Appendix (VI),
 we perform a numerical analysis and provide supporting empirical evidence that is consistent with the
formal theorems that were key to establishing the correctness of the RH. In
particular, the numerical results of the statistical tests performed here are firmly
consistent with the proposition (formally proved in the paper cited above)
that the values taken on by the Liouville function over large sequences of
consecutive integers are random. By performing this exhaustive numerical analysis
  and statistical study we feel that we have  provided a clearer understanding of the Riemann Hypothesis and its proof.
  
\begin{center}
{\Large 1. Introduction}
\end{center}

The Riemann zeta function, $\zeta (s)$, is defined by 
\begin{equation}
\zeta (s)=\sum\limits_{n=1}^{\infty }\frac{1}{n^{s}},
\end{equation}

\noindent where $n$ is a positive integer and $s$ is a complex variable,
with the series being convergent for $Re(s)>1$. This function has zeros
(referred to as the trivial zeros) at the negative even integers $%
-2,-4,\ldots $. It has been shown\footnote{%
This was first proved by Hardy (1914).} that there are an infinite number of
non-trivial zeros on the critical line at $Re(s)=1/2$. Riemann's Hypothesis
(RH), which has long remained unproven, claims that all the nontrivial zeros
of the zeta function lie on the critical line. The main paper contains the proof [MP and A's]

In this technical note, we provide a more concrete understanding and
appreciation of the steps involved in the proof of the Riemann Hypothesis by
supplying supporting empirical evidence for those various theorems which were proved and which had played
 a key role in the proof of the RH. In what follows we first
give a brief summary of how the RH was proved in [MP and A's].

The proof followed the primary idea that if the zeta function has zeros only
the critical line, then the function $F(s)$ $\equiv \zeta (2s)/\zeta (s)$
cannot be analytically continued to the left from the region $Re(s)>1$,
where it is analytic, to the left of $Re(s)<1/2$. This point was recognized
by Littlewood as far back as 1912.\footnote{It may be noted that Littlewood studied the function
 $ 1/\zeta(s)$ whereas we, in our analysis study $F(s)$ $\equiv \zeta (2s)/\zeta (s)$. This has made things simpler.}   The function $F(s)$ can be expressed as
(see Titchmarsh (1951, Ch. 1)):

\noindent 
\begin{equation}
F(s)=\sum\limits _{n=1}^{\infty}\frac{\lambda(n)}{n^{s}},
\end{equation}

\noindent where $\lambda (n)$ is the Liouville function defined by $\lambda
(n)=(-1)^{\omega (n)}$, with $\omega (n)$ being the total number of prime
numbers in the factorization of $n$, including the multiplicity of the
primes. The proof of RH in [MP and A's] requires also the summatory Liouville
function, $L(N)$, which is defined as: 
\begin{equation}
L(N)=\sum\limits_{n=1}^{N}\,\lambda (n)
\end{equation}%

 The proof crucially depends on showing that the function $F(s) = \zeta(2s) / \zeta(s)$, has poles only on the critical line $s = 1/2+i \sigma$, which translates to zeros of $\zeta(s)$, on the self same critical line $s = 1/2 + i \sigma$, because all the values of $s$ which appear as poles of $F(s)$ are actually zeros of $\zeta(s)$, except for $s=1/2$. Since, the trivial zeros of $\zeta(s)$ which occur at $s = -2,-4,-6.... $ that is negative even integers, conveniently cancel out from numerator and denominator of the expression in $F(s)$), leaving only the non trivial zeros, also the pole of $\zeta(2s)$ will appear as a pole of $F(s)$, at $s = 1/2$.  So it just remains to show that all the poles of $F(s)$ lie on the critical line. This was  the Primary task of the paper.

The crucial condition then is that $F(s)$ is not continuable to the left of $%
Re(s)<1/2$, and therefore that the zeta function have zeros only on the
critical line,\footnote{Riemann had already shown that symmetry conditions ensure that there will be no 
zeros $ 0 < Re(s)< 1/2$ if it is found that there are no zeros in the region $1/2 < Re(s) < 1$ }  is that the asymptotic limit of the summatory Liouville
function be $|L(N)|\sim C\,N^{1/2}$, where $C$ is a constant. Therefore, to
provide a rigorous proof of the validity of the Riemann Hypothesis, [MP and A's] 
investigated the asymptotic limit of $L(N)$. The work involved the
establishment of several relevant theorems, which were then invoked to
eventually prove the RH to be correct.

We now state some of these important theorems\footnote{In addition to the theorems given below, a necessary  theorem which states that: The sequence $\lambda(1),\lambda(2),\lambda(3),...,\lambda(n),...,$
is asymptotically non-cyclic, (i.e. it will never repeat), was also proved, in [MP and A's], the theorems are numbered differently}).

\ 

\ \textbf{\textit{Theorem}}\textbf{\ 1:}

In the set of all positive integers, for every integer which has an even
number of primes in its factorization there is another unique integer (its
twin) which has an odd number of primes in its factorization.

\underline{Remark}: Theorem 1 gives us the formal result that $Pr(\lambda
(n)=+1)=Pr(\lambda (n)=-1)=1/2$ , where Pr denotes probability. That is, the 
$\lambda $-function behaves like an `ideal coin'.

\ 

\ \textbf{\textit{Theorem}}\textbf{\ 2:}

Consider the sequence $S_{-}(N)$ comprising $\mu (N)$ consecutive positive
integers, defined by

\noindent $S_{-}(N)=\{N-\mu (N)+1,N-\mu (N)+2,N-\mu (N)+3,...,N\}$,
 \noindent where $\mu (N)=\sqrt{N}$. Then every number in $S_{-}(N)$\ will
differ in its prime factorization from that of every other 
number in $S_{-}(N)$ by at least one 
distinct prime.\footnote{For example, if two numbers $c$ and $d$ in $S$ are 
factorized as $c=p_{1}^{e_{1}}p_{2}^{e_{2}}$ and $d=p_{3}^{e_{3}}p_{4}^{e_{4}}$ then at
least one of the primes $p_{3}$ or $p_{4}$ will be different from $p_{1}$ or 
$p_{2}$.}

\underline{Remark}: It is not hard to prove that the 
sequence $S_{+}(N)\equiv \{N+1,N+2,N+3,........,N+\mu (N)\}$ of length $\mu (N)$ also
behaves similarly. That is, every member of $S_{+}(N)$ differs from every
other member by at least one prime in its factorization. This implies, interestingly, 
that by choosing $N$ to be consecutive perfect squares, the entire set of positive integers 
can be envisaged as a union of mutually exclusive sequences like $S_{-}(N)$ and $S_{+}(N)$.

It follows that the 
$\lambda $-values in the sequences $S_{-}^{\lambda }(N)\equiv \{\lambda
(N-\mu (N)+1),\lambda (N-\mu (N)+2),...,\lambda (N)\}$ and $S_{+}^{\lambda
}(N)\equiv \{\lambda (N+1),\lambda (N+2),\lambda (N+3),...,\lambda (N+\mu
(N))\}$, in which each member has a value either $+1$ or $-1$, would also
appear randomly and be statistically similar to sequences of coin tosses.

Since the number of members in the sequences $S_{-}(N)$, $S_{-}^{\lambda }(N)$,
 $S_{+}^{\lambda }(N)$, and $S_{+}^{\lambda }(N)$ is given by 
$\mu (N)= \sqrt{N} \rightarrow \infty $ as $N \rightarrow \infty $, the behavior 
of the $\lambda $-values of very large integers should coincide with that of a
sequence of coin tosses. This intuition was formally confirmed in 
Appendix $V$.

\textbf{\textit{Theorem}}\textbf{\ 3:}

The summatory Liouville function takes the asymptotic form
$|L(N)|=C\,N^{1/2} $,$C$ is 
a constant. It can be shown that $C=\sqrt{\frac{2}{\pi }}$. 
It may be mentioned here that Littlewood's condition is fairly tolerant: As long as asymptotically, for  large $N$,  $|L(N)| = C\,N^{1/2} $, and C is any finite constant, R.H. follows. This `tolerance' is reflected in the value of $\chi^2$ (below) as may be deduced,  after a study of the following.

\
\underline{Remark}: The form of the summatory Liouville function in
Theorem 3 is precisely what we would expect for a sequence of unbiased coin
tosses. This, along with a sufficient condition derived by Littlewood
(1912), shows that $F(s)$ is analytic for $Re(s)>1/2$ and $Re(s)<1/2$,
thereby leaving the only possibility that the non-trivial zeros of $\zeta (s)
$ can occur only on the critical line $Re(s)=1/2$.

\ In the following sections, by comparing the $\lambda $-sequences obtained
for large sets of consecutive integers with (binomial) sequences of coin
tosses, we show that the statistical distributions of the two sets of
sequences are consistent with the claims of the above theorems. To this end,
we apply Pearson's `Goodness of Fit' $\chi ^{2}$ test. The software program 
\textit{Mathematica}\textbf{\ }developed by Wolfram has been used in this
technical report to aid in the prime factorization of the large numbers that
this exercise entails. 

The compelling bottom line that emerges from this empirical study is that it
is extremely unlikely, in fact statistically impossible, that for large $N$, 
the sequences of $\lambda $%
-values can differ from sequences of coin tosses. It is this behavior of the
Liouville function, recall, that delivers Theorem 3 above. And this Theorem,
in turn, nails down all the non-trivial zeros of the zeta function to the
critical line [MP and A's].

\begin{center}
{\Large 2.  $\chi ^{2}$ Fit of a $\lambda $-Sequence}
\end{center}

In this section we will derive an expression of how closely a $\lambda $
sequence corresponds to a binomial sequence (coin tosses). We follow the
exposition given in Knuth (1968, Vol. 2, Ch. 3); and then derive a very important
expression for a $\chi ^{2}$ fit of a $\lambda $-Sequence, given by Eq.(14.9) below.

Suppose we are given a sequence, $T(N_{0},N)$, of $N$ consecutive integers
starting from $N_{0}$:

$T(N_{0},N)=\{N_{0},\,N_{0}+1,\,N_{0}+2,\,N_{0}+3,.......,\,N_{0}+N-1\}$

\noindent and the sequence, $\Lambda (N_{0},N)$, of the corresponding $%
\lambda $-values:

$\Lambda (N_{0},N)=\{\lambda (N_{0}),\,\lambda (N_{0}+1),\,$ $\lambda
(N_{0}+2),\,\lambda (N_{0}+3),..$ $,\,\lambda (N_{0}+N-1)\}$.

We ask how close in a statistical sense the sequence $\Lambda (N_{0},N)$ is
to a sequence of coin tosses or, in other words, a binomial sequence. By
identifying $\lambda (n) =1$ as $Head$ and $\lambda (n) =-1$ as $Tail$, for the $n^{th}$ `toss', we may
perform this comparison. If this is really the case then statistically $%
\Lambda (N_{0},N)$ should resemble a binomial distribution, we can then
compute the $\chi ^{2}$ statistic as follows.

\begin{equation}
\chi ^{2}(N)=\frac{(P-E_{P})^{2}}{E_{P}}+\frac{(M-E_{M})^{2}}{E_{M}},
\end{equation}

\noindent where $P$ and $M$ are the actual number of $+1$s ($Heads$) and $-1$%
s ($Tails$), respectively, in the $\Lambda (N_{0},N)$ sequence, $E_{P}$ and $%
E_{M}$ are the expectations of the number of $+1$s and $-1$s in the
probabilistic sense. From Theorem 1 it immediately follows that, for large $N
$,

\begin{equation}
E_{P}=E_{M}=N/2.
\end{equation}

We define $L(N_{0},N)$ as the additional contribution to the summatory
Liouville function of $N$ consecutive integers starting from $N_{0}$:

\begin{equation}
L(N_{0},N)=\sum_{n=N_{0}}^{N_{0}+N-1}\lambda (n).
\end{equation}

For brevity we will denote $\widehat{L}\equiv L(N_{0},N)$ and since (6)
contains $P$ terms which are equal to $+1$s and M terms which are equal to $%
-1$s, we can write:

\begin{equation}
P-M=\widehat{L},
\end{equation}

and

\begin{equation}
P+M=N.
\end{equation}

Using (7)and (8) we see that $P=(N+\widehat{L})/2$ and $M=(N-\widehat{L})/2$
and from (5) we deduce $P-E_{P}=\widehat{L}/2$ and $M-E_{M}=-\widehat{L}/2$
and thus equation (4) gives us the very important $\chi ^{2}$ relation which
is satisfied by every $\Lambda (N_{0},N)$ sequence involving the
factorization of $N$ consecutive integers starting from $N_{0}$:

\begin{equation}
\chi ^{2}(N)=\frac{[L(N_{0},N)]^{2}}{N}.
\end{equation}

Note that the it was possible to derive an expression for $\chi ^{2}$ for
large$N$ only because of Theorems 1, 2, and 3. Now we particularly choose $N$
to be the square of an integer and the sequence of length $\mu (N)=\sqrt{N}$
starting from the integer $N_{0}=N-\sqrt{N}+1$ and then taking the $\sqrt{N}$
consecutive terms of the $\lambda $-sequence, we obtain

\begin{equation}
\Lambda (N_{0},\sqrt{N})=\{\lambda (N-\sqrt{N}+1),\,\lambda (N-\sqrt{N}%
+2),\,\lambda (N-\sqrt{N}+3),\,...,\lambda (N)\,\},
\end{equation}

\noindent and the corresponding $\chi ^{2}(\sqrt{N})$ for such a sequence
(which is of length $\sqrt{N})$ can be obtained from Theorem 3 and (9) as%
\begin{eqnarray}
\chi ^{2}(\sqrt{N}) &=&\frac{[C\,\sqrt{\sqrt{N}}]^{2}}{\sqrt{N}}  \nonumber
\\
&=&C^{2}.
\end{eqnarray}

Equation (11) of course, should be interpreted as the average value of a
sequence such as $\Lambda (N_{0},\sqrt{N})$ of length $\sqrt{N}$ given in the
expression (10). In this report we perform the $\chi ^{2}$ `Goodness of Fit'
tests for very many sequences of the type $\Lambda (N_{0},\sqrt{N})$ with
varying lengths and very large values of $N$ to examine whether these
sequences are statistically indistinguishable from coin tosses. In this
manner, we provide empirical support for the claims of the theorems formally
proved in [MP and A's] and, therefore, for the proof of the Riemann Hypothesis.

$ $

{\Large 3. Numerical Analysis of Sequence $\Lambda (N_{0},\protect\sqrt{N})$ and
its $\protect\chi ^{2}$ Fit}

$ $
In this section, we consider sequences of length $\sqrt{N}$, starting from $N_{0}=N-\sqrt{N}+1$
or $N+1$ where $N$ is
a perfect square. We use \textit{Mathematica} to compute $L(N_{0},N)$.%
\footnote{%
A typical Mathematica command which calculates the expression
 $\sum_{n=J}^{K}\lambda (n)$ is:
 
 Plus{[}LiouvilleLambda{[}Range{[}J,K{]}{]}{].}For instance, the command which sums 
 
the $\lambda (n)$ from $n=25,000,001$ to $25,005,000$ is: 

Plus{[}LiouvilleLambda{[}Range{[}$25000001$, $25005000${]}{]}{], }
which will give the answer = $-42$.}

In the table below we list the sequences in the following format. We define
the sequences:

\begin{equation}
S_{-}(N)=\{N-\sqrt{N}+1,\,N-\sqrt{N}+2,\,...,N\},
\end{equation}

\begin{equation}
S_{+}(N)=\{N+1,\,N+2,\,...,\,N+\sqrt{N}\},
\end{equation}

\noindent and the partial sums of the $\lambda $s of the two sequences
defined above are defined by the expressions: 
\begin{equation}
L(S_{-})\equiv L(N-\sqrt{N}+1,N)=\lambda (N-\sqrt{N}+1)+\lambda (N-%
\sqrt{N}+2)+...+\lambda (N),
\end{equation}

\begin{equation}
L(S_{+})\equiv L(N+1,N+\sqrt{N})=\lambda (N+1)+\lambda (N+2)+...+\lambda (N+%
\sqrt{N}).
\end{equation}

The formal proof of the Riemann Hypothesis in [MP and A's] proceeded as follows. The
sequences $\Lambda (N-\sqrt{N}+1,\sqrt{N})$ and $\Lambda (N+1,\sqrt{N})$ 
were shown to behave like coin tosses for
every $N$ (large) over sequences of length $\sqrt{N}$, where $N$ is taken to
be a perfect square. On taking $N$ to be consecutive perfect squares, the
lengths of the consecutive sequences naturally increase. Using this
procedure, we obtain sequences that can span the \textit{entire set} of
positive integers (consult the first five columns of Tables 1.1 to 1.4).
Since the $\lambda $s within each segment behave like coin tosses, from the
work of Chandrashekar (1943) it follows that the summatory Liouville
function $L(N)$ must behave like $C\sqrt{N}$ as $N\rightarrow \infty $. The
validity of RH follows, by Littlewood's Theorem, from the fact that $F(s)$ cannot then be continued to
the left of the critical line $Re(s)=1/2$ because of the appearance of poles
in $F(s)$ on the line, each pole corresponding to a zero of the zeta
function $\zeta (s)$.

\ 

\textbf{Statistical Tests}

We shall now test the following null hypothesis $H_{0}$ against the
alternative hypothesis $H_{1}$ in the following generic forms:

$H_{0}$: The sequence $\Lambda (N_{0},N)$ has the same statistical
distribution as a corresponding sequence of coin tosses (i.e. binomial
distribution with $Prob(H)=Prob(T)=1/2$).

$H_{1}$: The sequence $\Lambda (N_{0},N)$ has a different statistical
distribution than a corresponding sequence of coin tosses (i.e. binomial
distribution with $Prob(H)=Prob(T)=1/2$).

The critical value for chi square is $\chi _{crit}^{2}=3.84$, for the
standard $0.05$ level of significance. (In our case, the relevant degrees of
freedom equal to $1$.) Assuming that $H_{0}$ is true, if chi square is less
than $\chi _{crit}^{2}$ the null hypothesis is accepted.

\ 

It should be noted that the tests conducted here are not merely exploratory
statistical exercises to discern possible patterns in the $\lambda $%
-sequences. Rather, the tests here are informed by theory. We have formally
shown in [MP and A's] that, over the set of positive integers, the probability that $%
\lambda $ takes on the value $+1$ or $-1$ with equal probability and that,
over sequences that are increasing in $N$, the $\lambda $ draws are random.
Thus statistical evidence consistent with these claims merely bolster what
has already been formally demonstrated.

The behavior of the $\Lambda (N_{0},\sqrt{N})$ sequences are verified to be
indeed like coin tosses for a very large number of cases and the results are
summarized in the tables below. Let us take an example from Table 1.1. The
third row gives the $\chi ^{2}$ result for the sequence of length $1001$,
starting from $1001001$. We can factorize each of these numbers as:

$1001001=3\,\times 333667;\,1001002=2\,\times \,500501;\,1001003=prime;$

$1001004=2^{2}\times \,3\,\times \,83417;$..........$,1001999=41\,\times \,24439;$

$\,1002000=2^{4}\times \,3\,\times \,5^{3}\times.\,167;$ $\,1002001=7^{2}\times \,11^{2}\times \,13^{2}$ 

and hence we can
evaluate the corresponding $\lambda $-sequence, by using the definition $%
\lambda (n)=(-1)^{\omega (n)}$, with $\omega (n)$ being the total number of
prime numbers (multilpicities included) in the factorization of $n$. We find
that:

$\Lambda (1001001,1001)=\{\lambda (1001001),$$\,$$\lambda (1001002),$$%
\,...,\,\lambda (1002000),$$\,$$\lambda (1002001)\}$

$=\{1,1,-1,1,1,1,1,-1,1,-1,1,................,1,-1,1\}.$

The partial sum of all the $1001$ $\lambda $s shown in the sequence above
adds up to $49$. We then estimate how close the sequence $\Lambda
(1001001,1001)$ is to a Binomial distribution, i.e. of $1001$ consecutive
coin tosses. The observed value $\chi ^{2}=2.4$ for this sequence of $%
\lambda $s is well below the critical value $\chi _{crit}^{2}=3.84$ (for a
one degree of freedom) at the standard significance level of $0.05$. Thus
the sequence $\Lambda (1001001,1001)$ is statistically indistinguishable
from a Binomial distribution obtained by $1001$ consecutive coin tosses if
we consider $Head=$ $+1$ and $Tail=-1$. In fact, it so happens that out of
the $10$ sequences shown in Table 1.1 this chosen example has the largest
value of $\chi ^{2}$; the other sequences have a much lower $\chi ^{2}$
value and the average value is $0.653$ which hovers around the predicted
average $C^{2}=\frac{2}{\pi }=0.637$. We see that the null hypothesis would
be accepted even if the significance level were at $0.10$, for which $\chi
_{crit}^{2}=2.71$.

We have calculated the $\chi ^{2}$ for larger and larger sequences see
Tables1.2, Tables 1.3 and Tables 1.4\textbf{\ }for even very large numbers $%
\sim 10^{10}$ and sequences involving $10^{5}$consecutive integers in each
case the sequences $\Lambda (N_{0},\sqrt{N})$ behave like coin tosses thus
lending emphatic empirical support consistent with the Theorems proved in [MP and A's],.

\medskip{}
\textbf{TABLE 1.1 Sequence of Consecutive Integers of Type $S_{-}(N)$
and $S_{+}(N)$ of Length 1000}

\medskip{}

\begin{tabular}{|c|c|c|c|c|c|c|}
\hline 
No & $Type$ of S(N) & $\sqrt{N}$ & From & to  & L(S) & $\chi^{2}$\tabularnewline
\hline 
\hline 
1. & S\_ & 1000 & 999,001 & 1,000,000 & 6 & 0.036\tabularnewline
\hline 
2. & $S_{+}$ & 1000 & 1,000,001 & 1,001,000 & 10 & 0.100\tabularnewline
\hline 
3. & S\_ & 1001 & 1,001,001 & 1,002,001 & 49 & 2.400\tabularnewline
\hline 
4. & $S_{+}$ & 1001 & 1,002,002 & 1,003,002 & -37 & 1.368\tabularnewline
\hline 
5. & S\_ & 1002 & 1,003,003 & 1,004,004 & -12 & 0.144\tabularnewline
\hline 
6. & $S_{+}$ & 1002 & 1,004,005 & 1,005,006 & -28 & 0.780\tabularnewline
\hline 
7. & S\_ & 1003 & 1,005,007 & 1,006,009 & 3 & 0.009\tabularnewline
\hline 
8. & $S_{+}$ & 1003 & 1,006,010 & 1,007,012 & -39 & 1.516\tabularnewline
\hline 
9. & S\_ & 1004 & 1,007,013 & 1,008,016 & 12 & 0.143\tabularnewline
\hline 
10. & $S_{+}$ & 1004 & 1,008,017 & 1,009,020 & 6 & 0.036\tabularnewline
\hline 
 &  &  &  &  &  & \tabularnewline
\hline 
 & \textbf{MEAN} & \textbf{$\chi^{2}$ FROM} & \textbf{999,001} & \textbf{1,009,020} & \textbf{=} & \textbf{0.653}\tabularnewline
\hline 
\end{tabular}

\newpage

\textbf{TABLE 1.2 Sequence of Consecutive Integers of Type $S_{-}(N)$
and $S_{+}(N)$ of Length 5000}

\medskip{}

\begin{tabular}{|c|c|c|c|c|c|c|}
\hline 
No & $Type$ of S(N) & $\sqrt{N}$ & From & to  & L(S) & $\chi^{2}$\tabularnewline
\hline 
\hline 
1. & S\_ & 5000 & 24,995,001 & 25,000,000 & 0 & 0.0\tabularnewline
\hline 
2. & $S_{+}$ & 5000 & 25,000,001 & 25,005,000 & -42 & 0.353\tabularnewline
\hline 
3. & S\_ & 5001 & 25,005,001 & 25,010,001 & -27 & 0.148\tabularnewline
\hline 
4. & $S_{+}$ & 5001 & 25,010,002 & 25,015,002 & -103 & 2.12\tabularnewline
\hline 
5. & S\_ & 5002 & 25,015,003 & 25,020,004 & -76 & 1.155\tabularnewline
\hline 
6. & $S_{+}$ & 5002 & 25,020,005 & 25,025,006 & 48 & 0.461\tabularnewline
\hline 
7. & S\_ & 5003 & 25,025,007 & 25,030,009 & -13 & 0.034\tabularnewline
\hline 
8. & $S_{+}$ & 5003 & 25,030,010 & 25,035,012 & 119 & 2.831\tabularnewline
\hline 
9. & S\_ & 5004 & 25,035,013 & 25,040,016 & 124 & 3.072\tabularnewline
\hline 
10. & $S_{+}$ & 5004 & 25,040,017 & 25,045,020 & 62 & 0.768\tabularnewline
\hline 
 &  &  &  &  &  & \tabularnewline
\hline 
 & \textbf{MEAN} & \textbf{$\chi^{2}$ FROM} & \textbf{24,995,001} & \textbf{25,045,020} & \textbf{=} & \textbf{1.094}\tabularnewline
\hline 
\end{tabular}

\medskip{}
\medskip{}

\textbf{TABLE 1.3 Sequence of Consecutive Integers of Type $S_{-}(N)$
and $S_{+}(N)$ of Length 10,000}

\medskip{}

\begin{tabular}{|c|c|c|c|c|c|c|}
\hline 
No & $Type$ & $\sqrt{N}$ & From & to  & L(S) & $\chi^{2}$\tabularnewline
\hline 
\hline 
1. & S\_ & 10000 & 99,990,001 & 100,000,000 & -146 & 2.132\tabularnewline
\hline 
2. & $S_{+}$ & 10000 & 100,000,001 & 100,010,000 & -88 & 0.774\tabularnewline
\hline 
3. & S\_ & 10001 & 100,010,001 & 100,020,001 & -11 & 0.012\tabularnewline
\hline 
4. & $S_{+}$ & 10001 & 100,020,002 & 100,030,002 & -43 & 0.185\tabularnewline
\hline 
5. & S\_ & 10002 & 100,030,003 & 100,040,004 & 8 & 0.064\tabularnewline
\hline 
6. & $S_{+}$ & 10002 & 100,040,005 & 100,050,006 & 36 & 0.130\tabularnewline
\hline 
7. & S\_ & 10003 & 100,050,007 & 100,060,009 & 23 & 0.053\tabularnewline
\hline 
8. & $S_{+}$ & 10003 & 100,060,010 & 100,070,012 & -49 & 0.240\tabularnewline
\hline 
9. & S\_ & 10004 & 100,070,013 & 100,080,016 & -20 & 0.040\tabularnewline
\hline 
10.  & $S_{+}$ & 10004 & 100,080,017 & 100,090,020 & 112 & 1.254\tabularnewline
\hline 
 &  &  &  &  &  & \tabularnewline
\hline 
 & \textbf{MEAN} & \textbf{$\chi^{2}$ FROM} & \textbf{99,990,001 TO} & \textbf{100,090,020} & \textbf{=} & \textbf{0.488}\tabularnewline
\hline 
\end{tabular}

\medskip{}

\newpage

\textbf{TABLE 1.4 Sequence of Consecutive Integers of Type $S_{-}(N)$
and $S_{+}(N)$ of Length 100,000}

\medskip{}

\begin{tabular}{|c|c|c|c|c|c|c|}
\hline 
No & $Type$ & $\sqrt{N}$ & From & to  & L(S) & $\chi^{2}$\tabularnewline
\hline 
\hline 
1. & S\_ & 100,000 & 9,999,900,001 & 10,000,000,000 & -232 & 0.538\tabularnewline
\hline 
2. & S$_{+}$ & 100,000 & 10,000,000,001 & 10,000,100,000 & 340 & 1.15\tabularnewline
\hline 
3. & S\_ & 100,001 & 10,000,100,001 & 10,000,200,001 & -249 & 0.620\tabularnewline
\hline 
4. & S$_{+}$ & 100,001 & 10,000,400,005 & 10,000,500,006 & -115 & 0.132\tabularnewline
\hline 
5. & S\_ & 100,002 & 10,000,300,003 & 10,000,400,004 & 216 & 0.467\tabularnewline
\hline 
6. & S$_{+}$ & 100,002 & 10,000,400,005 & 10,000,500,006 & 456 & 2.08\tabularnewline
\hline 
7. & S\_ & 100,003 & 10,000,500,007 & 10,000,600,009 & -255 & 0.650\tabularnewline
\hline 
8. & S$_{+}$ & 100,003 & 10,000,600,010 & 10,000,700,012 & -235 & 0.552\tabularnewline
\hline 
9. & S\_ & 100,004 & 10,000,700,013 & 10,000,800,016 & -44 & 0.0194\tabularnewline
\hline 
10.  & S$_{+}$ & 100,004 & 10,000,800,017 & 10,000,900,020 & 202 & 0.408\tabularnewline
\hline 
11. & S\_ & 100,005 & 10,000,900,021 & 10,001,000,025 & -191 & 0.364\tabularnewline
\hline 
12. & S$_{+}$ & 100,005 & 10,001,000,026 & 10,001,100,030 & 475 & 2.26\tabularnewline
\hline 
13. & S\_ & 100,006 & 10,001,100,031 & 10,001,200,036 & 134 & 0.179\tabularnewline
\hline 
14. & S$_{+}$ & 100,006 & 10,001,200,037 & 10,001,300,042 & -66 & 0.0436\tabularnewline
\hline 
15. & S\_ & 100,007 & 10,001,300,043 & 10,001,400,049 & 427 & 1.82\tabularnewline
\hline 
16. & S$_{+}$ & 100,007 & 10,001,400,050 & 10,001,500,056 & -303 & 0.918\tabularnewline
\hline 
17. & S\_ & 100,008 & 10,001,500,057 & 10,001,600,064 & 276 & 0.762\tabularnewline
\hline 
18. & S$_{+}$ & 100,008 & 10,001,600,065 & 10,001,700,072 & -210 & 0.441\tabularnewline
\hline 
19. & S\_ & 100,009 & 10,001,700,073 & 10,001,800,081 & 267 & 0.713\tabularnewline
\hline 
20. & S$_{+}$ & 100,009 & 10,001,800,082 & 10,001,900,090 & 291 & 0.847\tabularnewline
\hline 
 &  &  &  &  &  & \tabularnewline
\hline 
 & \textbf{MEAN} & \textbf{$\chi^{2}$ FROM} & \textbf{9,999,900,001 TO} & \textbf{10,001,900,090} & \textbf{=} & \textbf{0.768}\tabularnewline
\hline 
\end{tabular}

\medskip{}

\medskip{}

{\Large 3.1 Sequences of Fixed Length Arbitrarily Positioned}

$ $

In this section we consider various segments of consecutive integers of a
fixed length but starting from an arbitrary integer. Even here we see that
the $\lambda $s within each segment behave like coin tosses and have the
same statistical properties.

We now calculate the $\chi ^{2}$ values of $\lambda $-sequences for a
sequence $S_{A}$ of consecutive integers, starting from an arbitrary number $%
N_{0}$ but all of a fixed length $M$:

\begin{equation}
S_{A}(N)=\{N_{0},\,N_{0}+1,\,N+2,\,N+3,...,\,N_{0}+M-1\}
\end{equation}

and

\begin{equation}
L(S_{A})=\lambda (N_{0})+\lambda (N_{0}+1)+\lambda (N_{0}+2)+\lambda
(N_{0}+3)+...+\lambda (N_{0}+M-1)
\end{equation}

\medskip {}The results, which are summarized in Table 2.1, again show that
the $\lambda $-sequences are statistically like coin tosses.

\textbf{TABLE 2.1 Sequence of Consecutive Integers of Type $S_{A}(N)$
and  of Length M = 1000}

\medskip{}

\begin{tabular}{|c|c|c|c|c|c|c|}
\hline 
No & $Type$ & M  & From $N_{0}$ & to $N_{0}+M-1$ & $L(S_{A})$ & $\chi^{2}$\tabularnewline
\hline 
\hline 
1. & $S_{A}$ & 1000 & 10,000,001 & 10,001,000 & 36 & 1.296\tabularnewline
\hline 
2. & $S_{A}$ & 1000 & 12,000,001 & 12,001,000 & 28 & 0.784\tabularnewline
\hline 
3. & $S_{A}$ & 1000 & 13.000,001 & 13,001,000 & -14 & 0.196\tabularnewline
\hline 
4. & $S_{A}$ & 1000 & 15,000,001 & 15,001,000 & 10 & 0.10\tabularnewline
\hline 
5. & $S_{A}$ & 1000 & 45,000,001 & 45,001,000 & -18 & 0.324\tabularnewline
\hline 
6. & $S_{A}$ & 1000 & 47,000,001 & 47,001,000 & -36 & 1.296\tabularnewline
\hline 
7. & $S_{A}$ & 1000 & 56,000,001 & 56,001,000 & 24 & 0.576\tabularnewline
\hline 
8. & $S_{A}$ & 1000 & 70,000,001 & 70,001,000 & -44 & 1.936\tabularnewline
\hline 
9. & $S_{A}$ & 1000 & 90,000,001 & 90,001,000 & 14 & 0.196\tabularnewline
\hline 
10. & $S_{A}$ & 1000 & 95,600,001 & 95,601,000 & 28 & 0.784\tabularnewline
\hline 
11. & $S_{A}$ & 1000 & 147,000,001 & 147,001,000 & -26 & 0.676\tabularnewline
\hline 
12. & $S_{A}$ & 1000 & 237,000,001 & 237,001,000 & -24 & 0.576\tabularnewline
\hline 
13. & $S_{A}$ & 1000 & 400,000,001 & 400,001,000 & 26 & 0.676\tabularnewline
\hline 
14. & $S_{A}$ & 1000 & 413,000,001 & 413,001,000 & 10 & 0.10\tabularnewline
\hline 
15. & $S_{A}$ & 1000 & 517,000,001 & 517,001,000 & 14 & 0.196\tabularnewline
\hline 
16. & $S_{A}$ & 1000 & 530,000,001 & 530,001,000 & -32 & 1.024\tabularnewline
\hline 
17. & $S_{A}$ & 1000 & 731.000,001 & 731,001,000 & 50 & 2.500\tabularnewline
\hline 
18. & $S_{A}$ & 1000 & 871,000,001 & 871,001,000 & -42 & 1.764\tabularnewline
\hline 
19. & $S_{A}$ & 1000 & 979,000,001 & 979,001,000 & -20 & 0.400\tabularnewline
\hline 
20. & $S_{A}$ & 1000 & 997,000,001 & 997,001,000 & 14 & 0.196\tabularnewline
\hline 
 &  &  &  &  &  & \tabularnewline
\hline 
 &  & \textbf{MEAN} & \textbf{$\chi^{2}$ OF ABOVE} & \textbf{20 SEGMENTS} & \textbf{=} & \textbf{0.780}\tabularnewline
\hline 
\end{tabular}

\medskip{}

\medskip{}

{\Large 3.2 Entire Sequences from $n=1$ to $n=N$, N large and calculation of $\chi^2$ for such sequences from $L(N)$}

$ $ 

It has been empirically verified in the literature that the summatory
Liouville Function $L(N)=\sum_{n=1}^{N}\lambda (n)$ fluctuates from positive
to negative values as $N$ increases without bound. We now investigate the $\chi ^{2}$ values
for such sequences,and use Eq.(9), so that we may see how these sequences behave like coin 
tosses.

In the Table 3.1 we use the values of $L(N)$ for various large values of $N$, which were found by
Tanaka (1980), the results depicted below reveal that the lambda sequences are statistically 
indistinguishable from the sequences of coin tosses over such large ranges of N from 1 to one billion.

\medskip{}

In the above we calculated $L(N)$ for various valies of $N$, however, 
 if we choose a value $N$ at which $L(N)$ is a local maximum or a local minimum then we would be
examining potential worst case scenarios for deviations of the $\lambda $s
from coin tosses because these are the values of $N$ that are likely to
yield the highest values of $\chi ^{2}$ (see equation (9)). It is
interesting to investigate if even for these special values of $N$ whether
the $\chi ^{2}$ is less than the critical value; if so, we would again have
statistical assurance that the entire sequence of $\lambda $s from $%
n=1,2,3,...$ behave like coin tosses.

We therefore use the $58$ largest values of $L(N)$ and the associated values
of $N$ reported in the literature by Borwein, Ferguson and Mossinghoff (2008),
 and perform our statistical exercise. See Table 3.2. We see that even for these
\textquotedblleft worst case scenario\textquotedblright\ values of $N$ 
the lambda sequences are statistically indistinguishable from the sequences
of coin tosses.

\textbf{TABLE 3.1 Values of L(N) at various large values of N }

(The values for N and L(N) are from Tanaka (1980))

\medskip{}

\begin{tabular}{|c|c|c|c|c|}
\hline 
No. & N & $L(N)=\sum_{n=1}^{N}\lambda(n)$ & $\chi^{2}$ & \tabularnewline
\hline 
1	&		100,000,000     &  -3884	       &   0.1508 \tabularnewline
2	&		200,000,000     &  -11126          &   0.6189 \tabularnewline
3	&		300,000,000     &  -16648         &   0.9238 \tabularnewline
4	&		400,000,000     &  -11200          &   0.3136 \tabularnewline
5	&		500,000,000     &  -18804          &   0.7072 \tabularnewline
6	&		600,000,000     &  -15350          &   0.3927 \tabularnewline
7	&		700,000,000     &  -25384          &   0.9204 \tabularnewline
8	&		800,000,000     &  -19292          &   0.4652 \tabularnewline
9	&		900,000,000     &  -4630           &   0.0238 \tabularnewline
10	&	  1,000,000,000     &  -25216          &   0.6358 \tabularnewline
	&		                 &					&			 \tabularnewline
	&	\textbf{MEAN}    & \textbf{$\chi^{2}$ OF ABOVE} =   &  0.5152                  \tabularnewline
	&		               &              &           \tabularnewline
\hline 
\end{tabular}$ $
\medskip{}

\newpage 

\textbf{TABLE 3.2 Values of L(N) at local Minima (Maxima) for very
Large N }

(The values for N and L(N) are from Borwein, Ferguson and Mossinghoff (2008))

\medskip{}

\begin{tabular}{|c|c|c|c|c|}
\hline 
No. & N & $L(N)=\sum_{n=1}^{N}\lambda(n)$ & $\chi^{2}$ & \tabularnewline
\hline 
1	&		293	            &  -21		       &   1.5051 \tabularnewline
2	&		468             &  -24             &   1.2308 \tabularnewline
3	&		684             &  -28             &   1.1462 \tabularnewline
4	&		1,132            &  -42             &   1.5583 \tabularnewline
5	&		1,760            &  -48             &   1.3091 \tabularnewline
6	&		2,804            &  -66             &   1.5535 \tabularnewline
7	&		4,528            &  -74             &   1.2094 \tabularnewline
8	&		7,027            &  -103            &   1.5097 \tabularnewline
9	&		9,840            &  -128             &   1.665 \tabularnewline
10	&		24,426           &  -186            &   1.4164 \tabularnewline
11	&		59,577           &  -307             &   1.582 \tabularnewline
12	&		96,862           &  -414            &   1.7695 \tabularnewline
13	&		386,434          &  -698            &   1.2608 \tabularnewline
14	&		614,155          &  -991            &   1.5991 \tabularnewline
15	&		925,985          &  -1,253           &   1.6955 \tabularnewline
16	&		2,110,931         &  -1,803             &   1.54 \tabularnewline
17	&		3,456,120         &  -2,254             &   1.47 \tabularnewline
18	&		5,306,119         &  -2,931            &   1.619 \tabularnewline
19	&		5,384,780         &  -2,932           &   1.5965 \tabularnewline
20	&		8,803,471         &  -3,461           &   1.3607 \tabularnewline
\hline 
\end{tabular}$ $

\medskip{}
\newpage

\textbf{TABLE 3.2 (Cont'd) Values of L(N) at local Minima (Maxima) for very
Large N }

(The values for N and L(N) are from Borwein, Ferguson and Mossinghoff (2008))

\medskip{}

\begin{tabular}{|c|c|c|c|c|}
\hline 
No. & N & $L(N)=\sum_{n=1}^{N}\lambda(n)$ & $\chi^{2}$ & \tabularnewline
\hline 

21	&		12,897,104        &  -4,878            &   1.845 \tabularnewline
22	&		76,015,169        &  -10,443          &   1.4347 \tabularnewline
23  &       184,699,341       &  -17,847          &   1.7245 \tabularnewline
24	&		281,876,941       &  -19,647          &   1.3694 \tabularnewline
25	&		456,877,629       &  -28,531          &   1.7817 \tabularnewline
26	&		712,638,284       &  -29,736          &   1.2408 \tabularnewline
27	&		1,122,289,008      &  -43,080          &   1.6537 \tabularnewline
28	&		1,806,141,032      &  -50,356          &   1.4039 \tabularnewline
29	&		2,719,280,841      &  -62,567          &   1.4396 \tabularnewline
30	&		3,847,002,655      &  -68,681          &   1.2262 \tabularnewline
31	&		     4,430,947,670  &  -73436          &   1.2171 \tabularnewline
32	&		     6,321,603,934  &  -96,460          &   1.4719 \tabularnewline
33	&		    10,097,286,319  &  -123,643         &   1.514  \tabularnewline
34	&		    15,511,912,966  &  -158,636         &   1.6223 \tabularnewline
35	&		    24,395,556,935  &  -172,987         &   1.2266 \tabularnewline
36	&		    39,769,975,545  &  -238,673         &   1.4324 \tabularnewline
37	&		    98,220,859,787  &  -365,305         &   1.3586 \tabularnewline
38	&		   149,093,624,694  &  -461,684         &   1.4296 \tabularnewline
39	&		   217,295,584,371  &  -598,109         &   1.6463 \tabularnewline
40	&		   341,058,604,701  &  -726,209         &   1.5463 \tabularnewline
41	&		   576,863,787,872  &  -900,668         &   1.4062 \tabularnewline
42	&		   835,018,639,060  &  -1,038,386        &   1.2913 \tabularnewline
43	&		 1,342,121,202,207  &  -1,369,777        &   1.398  \tabularnewline
44	&		 2,057,920,042,277  &  -1,767,635        &   1.5183 \tabularnewline
45	&		 2,147,203,463,859  &  -1,784,793        &   1.4836 \tabularnewline
46	&		 3,271,541,048,420  &  -2,206,930        &   1.4888 \tabularnewline
47	&		 4,686,763,744,950  &  -2,259,182        &   1.089  \tabularnewline
48	&		 5,191,024,637,118  &  -2,775,466        &   1.4839 \tabularnewline
49	&		 7,934,523,825,335  &  -3,003,875        &   1.1372 \tabularnewline
50	&		 8,196,557,476,890  &  -3,458,310        &   1.4591 \tabularnewline
51	&		12,078,577,080,679  &  -4,122,117        &   1.4068 \tabularnewline
52	&		18,790,887,277,234  &  -4,752,656        &   1.2021 \tabularnewline
53	&		20,999,693,845,505  &  -5,400,411        &   1.3888 \tabularnewline
54	&		29,254,665,607,331  &  -6,870,529        &   1.6136 \tabularnewline
55	&		48,136,689,451,475  &  -7,816,269        &   1.2692 \tabularnewline
56	&		72,204,113,780,255  &  -11,805,117       &   1.9301 \tabularnewline
57	&	   117,374,745,179,544  &  -14,496,306       &   1.7904 \tabularnewline
58	&	   176,064,978,093,269  &  -17,555,181       &   1.7504 \tabularnewline

\hline 
\end{tabular}$ $

\medskip{}
The empirical evidence provided here is very comprehensive: it examines the
statistical behavior of the Liouville function for large  segments 
of consecutive integers (e.g.Table 1.4). We have also considered the entire series
of $\lambda(n)$ from the values of $n=1$ to $n=N= 176$ trillion - as high as any available studies 
in the literature have gone. And yet, the $\lambda-$sequences consistently show themselves, in
rigorous statistical tests, to be indistinguishable from sequences of coin tosses, hence providing overwhelming statistical evidence in support of Littlewood's condition that as $N \rightarrow \infty$, $L(N) = C. \sqrt (N)$, (where C is finite)  and thus declaring that the non-trivial zeros of the zeta function, $\zeta(s)$, must all necessarily lie on the critical line $Re(s)= 1/2$.

{\Large 4. Concluding Note}

In this Appendix VI we have provided compelling, comprehensive numerical and statistical evidence that is consistent with the Theorems that were instrumental in the formal validation of the Riemann Hypothesis in [MP and A's].

It is hoped that a perusal of this section (Appendix 6)report offers some insight into, and understanding of, why the Riemann Hypothesis is correct. It should be noted that, while the results presented here are perfectly consistent with the theoretical results in [MP and A's],  they obviously do not prove (in a strict mathematical sense, because of the statistical nature of the study), the Riemann Hypothesis. For the formal proof, the rigorous mathematical analysis in the main paper needs to be consulted.

$ $

\section{References}
$ $
1. Apostol, T.,(1998) \textit{Introduction to Analytic Number Theory}%
, Chapter 2, pp. 37-38,  Springer International, Narosa Publishers, New Delhi.

2. Borwein, P., Choi,S., Rooney,B., and Weirathmueller,A., (2006) \textit{%
The Riemann Hypothesis}, Springer. Also see slides 28 to 32 in Lecture:

 http://www.cecm.sfu.ca/personal/pborwein/SLIDES/RH.pdf

3. Chandrasekhar S., (1943)'Stochastic Problems in Physics and Astronomy',Rev.of Modern Phys. 
vol 15, no 1, pp1-87 

4. Denjoy, A., (1931) L' Hypothese de Riemann sur la distribution des zeros.. C.R. Acad. Sci. Paris 192, 656-658

5. Edwards, H.M., (1974), \textit{Riemann's Zeta Function}, Academic
Press, New York.

6. Eswaran, K.,  (1990)`On the Solution of Dual Integral Equations Occurring in Diffraction Problems', Proc. of Roy. Soc. London, A 429, pp. 399-427

7. Hardy, G.H., (1914), `Sur les zeros de la fonction $\zeta (s)$ de
Riemann,' \textit{Comptes Rendus de l'Acad. des Sciences (Paris)}, 158,
1012--1014.

8. Knuth D.,(1968) `Art of Computer Programming', vol 2,Chap 3. (1968)Addison Wesley

9. A. Khinchine. ``Uber einen Satz der Wahrscheinlichkeitsrechnung'', Fundamenta Mathematicae 6 (1924): pp. 9-20

10. A. Kolmogoroff. ``Uber das Gesetz des iterierten Logarithmus''. Mathematische Annalen, 101: 126-135, 1929. (At the Göttinger DigitalisierungsZentrum web site)). Also see: https:\/\/en.wikipedia.org\/wiki\/Law of the iterated logarithm.

11. Landau,E., 1899 Ph. D Thesis, Friedrich Wilhelm Univ. Berlin 

12. Littlewood, J.E. (1912), `Quelques consequences
de l'hypoth`ese que la fonction $\zeta (s)$ de Riemann n'a pas de z%
eros dans le demi-plan $R(s)>1/2$,' \textit{Comptes Rendus
de l'Acad. des Sciences (Paris)}, 154, pp. 263--266. Transl. version,in Edwards op cit.

13. Riemann,B., (1892)in `Gessammeltz Werke', Liepzig, 1892, trans. `On the number of primes less than a given magnitude' available In Edwards' book, opcit pp299-305

14. Sautoy,Marcus du, (2003) 'Music of the primes' (Harper Collins)

15. Shannon, C.E., (1948)`A mathematical theory of communication',Bell Systems Technical Journal,
vol 27, pp 379-423, 623-656, July, October 1948. 

16.  Press,W.H, Teukolsky,S.A., Vetterling, W.T. and Flannery,B.P., Press etal (1986)`On Numerical Recipes in Fortran',Chap 7, Cambridge univ. Press, 1986

17. Whittaker, E.T. and Watson, G.N., (1989)`A Course of Modern Analysis', Universal Book Stall New Delhi 

18. Tanaka, Minoru, (1980), `A numerical investigation of the cumulative sum 
  of the Liouville function', Tokyo J. of Math. vol 3, No. 1, pp. 187-189.

19. Titchmarsh, E.C. (1951), \textit{The Theory of the Riemann Zeta-Function}%
, Clarendon Press, Oxford.

\end{document}